\documentclass[12pt]{article}
\usepackage{epsf,epsfig,amsmath,amsfonts,upgreek,xcolor}
\usepackage{graphicx}
\usepackage{multirow}
\usepackage{booktabs}
\usepackage{appendix}
\usepackage{mathrsfs}

\newtheorem{thm}{Theorem}

\newtheorem{lem}[thm]{Lemma}

\newtheorem{ass}{Assumption}
\newtheorem{de}{Definition}
\newtheorem{ex}{Example}
\newtheorem{rem}{Remark}

\newcommand{\be}{\begin{equation*}}
\newcommand{\ben}{\begin{equation}}
\newcommand{\ee}{\end{equation*}}
\newcommand{\een}{\end{equation}}

\newcommand{\bbe}{\operatorname{\mathbb{E}}}
\newcommand{\bbp}{\operatorname{\mathbb{P}}}
\newcommand{\bbr}{\mathbb{R}}

\newcommand{\acal}{\mathcal{A}}

\newcommand{\fcal}{\mathcal{F}}
\newcommand{\ical}{\mathcal{I}}

\newcommand{\lscr}{\mathscr{L}}

\newcommand{\di}{\mathrm{d}}

\newcommand{\xlb}{\underline{x}}
\newcommand{\xub}{\overline{x}}

\begin{document}

\title{The Solution to an Impulse Control Problem \\
Motivated by Optimal Harvesting}

\author{
{\sc Zhesheng Liu\footnote{
E-mail:~\texttt{z.liu64@lse.ac.uk}}
\ and
{\sc Mihail Zervos}\footnote{
E-mail:~\texttt{mihalis.zervos@gmail.com}}} \\
Department of Mathematics \\
London School of Economics \\
Houghton Street \\
London WC2A 2AE, UK}

\maketitle

\begin{abstract} \noindent
We consider a stochastic impulse control problem that is
motivated by applications such as the optimal
exploitation of a natural resource.
In particular, we consider a stochastic system whose
uncontrolled state dynamics are modelled by a non-explosive
positive linear diffusion.
The control that can be applied to this system takes the form of
one-sided impulsive action.
The objective of the control problem is to maximise a
discounted performance criterion that rewards the effect
of control action but involves a fixed cost at each time
of a control intervention.
We derive the complete solution to this problem
under general assumptions.
It turns out that the solution can take four qualitatively
different forms, several of which have not been
observed in the literature.
In two of the four cases, there exist only
$\varepsilon$-optimal control strategies.
We also show that the boundary classification of 0
may play a critical role in the solution of the problem.
Furthermore, we develop a way for establishing the
strong solution to a stochastic impulse control problem's
optimally controlled SDE.
\end{abstract}

%===============================================
\section{Introduction}

We consider a stochastic dynamical system whose controlled
state process is the strong solution to the SDE
\ben
\di X_t^\zeta = b(X_t^\zeta) \, \di t - \di \zeta_t  + \sigma
(X_t^\zeta) \, \di W_t , \quad X_{0-}^\zeta = x  > 0 , \label{SDE}
\een
where $W$ is a standard one-dimensional Brownian motion
and $\zeta$ is a controlled c\`{a}dl\`{a}g increasing piece-wise
constant process. 
The objective of the optimisation problem is to maximise
over all admissible processes $\zeta$ the performance criterion
\ben
J_x (\zeta) = \bbe_x \Biggl[ \int _0^\infty e^{-\Lambda_t^\zeta}
h (X_t^\zeta) \, \di t + \sum _{t \geq 0} e^{-\Lambda_t^\zeta}
\biggl( \int _0^{\Delta \zeta_t} k(X_{t-}^\zeta -u) \, \di u -
c {\bf 1} _{\{ \Delta \zeta_t > 0 \}} \biggr) \Biggr] ,  \label{JE}
\een
where $\Delta \zeta _t = \zeta _t - \zeta _{t-}$, with the
convention that $\zeta _{0-} = 0$, and
\ben
\Lambda _t^\zeta = \int _0^t r(X^\zeta_u) \, \di u .
\label{Lambda}
\een
Throughout the paper, we write $\bbe_x$ to denote
expectation so that we account for the dependence
of $X^{\zeta}$ on its initial value $x$.

Stochastic impulse control problems arise in various fields.
In the context of mathematical finance, economics and operations
research, notable contributions include
Harrison, Sellke and Tayor~\cite{HST83},
Harrison and Taksar~\cite{HT83},
Mundaca and {\O}ksendal~\cite{MO98},
Korn~\cite{K98, K99},
Bielecki and Pliska~\cite{BP00},
Cadenillas~\cite{C00},
Bar-Ilan, Sulem and Zanello~\cite{BSZ02},
Bar-Ilan, Perry and Stadje~\cite{BPS04},
Ohnishi and Tsujimura~\cite{OT06},
Cadenillas, Sarkar and Zapatero~\cite{CSZ07},
LyVath, Mnif and Pham~\cite{LVMP07},
and several references therein.
Also, impulse control models motivated by the optimal
management of a  natural resource have been studied
by
Alvarez~\cite{A041,A042},
Alvarez and Koskela~\cite{AK07} and
Alvarez and Lempa~\cite{AL08};
singular control versions of such models have
been studied by
Lungu and {\O}ksendal~\cite{LO01},
Framstad~\cite{F03},
Song, Stockbridge and Zhu~\cite{SSZ11},
Alvarez and Hening~\cite{AH22}
and several references therein.
In view of the wide range of applications, the general
mathematical theory of stochastic impulse control is
well-developed: apart from the contributions mentioned
above,
see also
Richard~\cite{R77},
Stettner~\cite{S83},
Lepeltier and Marchal~\cite{LM84},
Perthame~\cite{P84},
Egami~\cite{E08},
Davis, Guo and Wu~\cite{DGW10},
Djehiche, Hamad\`{e}ne and Hdhiri~\cite{DHH10},
Christensen~\cite{C14},
Helmes, Stockbridge and Zhu~\cite{HSZ15, HSZ24},
Menaldi and Robin~\cite{MR17},
Palczewski and Stettner~\cite{PS17},
Christensen and Strauch~\cite{CS23},
as well as the books by
Bensoussan and Lions~\cite{BL84},
Davis~\cite{D93},
{\O}ksendal and Sulem~\cite{OS07},
Pham~\cite{P09},
and several references therein.

The optimal management of a natural resource has
motivated the problem that we study here.
In this context, the state process $X^\zeta$ models the
population density of a harvested species, while
$\zeta_t$ is the cumulative amount of the species
that has been harvested by time $t$.
The constant $c>0$ models a fixed cost associated
with each harvesting cycle, while the function $k$
models the marginal profit arising from each harvest.
Furthermore, the function $h$ models the utility of the
harvested species, which could reflect the importance
of the species to its associated ecosystem.
Alternatively, the function $h$ can be used to model
the revenue or cost associated with running the
ecosystem.
Relative to related references, such as the ones
mentioned in the previous section, we generalise
by considering (a) state-space discounting,
(b) a state-dependent, rather than proportional, payoff
associated with each harvest size, and (c) a running
payoff such as the one modelled by the function
$h$.
On the other hand, the assumptions that we make
are of a rather similar nature.

In light of standard impulse control theory, a
``$\beta$-$\gamma$'' strategy should be a prime
candidate for an optimal one in the problem that we
study here.
Such a strategy is characterised by two points
$\gamma < \beta$ in $]0, \infty[$, which are both chosen
by the controller, and can be described
informally as follows.
If the state process  takes any value $x \geq \beta$,
then it is optimal for the controller to push it in an
impulsive way down to level $\gamma$.
On the other hand, the controller should wait and
take no action at all for as long as the state process
takes values in the interval $]0, \beta[$.

We show that a $\beta$-$\gamma$ strategy is indeed
optimal, provided that a critical parameter $\xlb$ is
finite and the fixed cost $c$ is sufficiently small (see
Case~I of Theorem~\ref{thm:SolPb-lem} in
Section~\ref{sec:solution}).
Otherwise, we show that only $\varepsilon$-optimal
strategies may exist (see Case~II or Case~IV of
Theorem~\ref{thm:SolPb-lem}) or that never making
an intervention may be optimal (see Case~III of
Theorem~\ref{thm:SolPb-lem}).
The absence of an optimal strategy in Case~IV of
Theorem~\ref{thm:SolPb-lem} is due to the relatively
rapid growth of the function $k$ at infinity.
It can therefore be eliminated if we make a suitable
additional growth assumption.
On the other hand, the absence of an optimal strategy
in Case~II of Theorem~\ref{thm:SolPb-lem} is due
to the nature of the problem that we solve.

The family of admissible controlled strategies that we
consider do not allow for the state process to hit the
boundary point 0 and be absorbed by it, which would
amount to ``switching off'' the system.
If we enlarged the set of admissible controls to allow
for such a possibility and 0 were a natural boundary point,
then we would face only the following difference:
a $\beta$-0 strategy would be optimal in Case~II of
Theorem~\ref{thm:SolPb-lem} and we would not need
to consider $\varepsilon$-optimal strategies.
On the other hand, the situation would be radically
different if 0 were an entrance boundary point:
in this case, $\beta$-0 strategies would become
an indispensable part of the optimal tactics.
We discuss these observations more precisely in
Remark~\ref{rem:nat/entr} at the end of
Section~\ref{sec:solution}.
To the best of our knowledge, this is the first
stochastic control problem in which the boundary
classification of the problem's state space has such
a fundamental influence on the problem's solution.
We do not investigate this issue any further because
this would require substantial extra analysis that
would go beyond the scope of the present article.

The evolution of an impulse control problem's state
process is quite intuitive, provided that the
corresponding uncontrolled dynamics are
well-posed.
For this reason, several references simply assume
the existence of such processes.
In the context of SDEs in $\bbr^d$, the state process
of an impulse control problem can be derived by
pasting together suitable strong solutions to the
underlying uncontrolled SDE with random initial
conditions (e.g., see Bensoussan and
Lions~\cite[Section~6.1.1]{BL84}).
In the context of general Markov processes,
the classical construction of an impulse control
strategy is substantially more technical and may
involve countable products of canonical spaces
(e.g., see Stettner~\cite{S83} and
Lepeltier and Marchal~\cite{LM84}).
If the uncontrolled state space process is a general
Markov process with continuous sample paths, then
comprehensive constructions of impulse control
models have been derived by Helmes, Stockbridge
and Zhu~\cite{HSZ24}.

Impulse control problems with SDEs in $\bbr^d$
can be formulated as in (\ref{SDE})--(\ref{Lambda}).
In itself such a formulation is straightforward.
Indeed, an SDE in $\bbr^d$ such as (\ref{SDE})
has a unique strong solution under suitable
Lipschitz assumptions on $b$ and $\sigma$ for
a wide class of controlled processes $\zeta$
(e.g., see Krylov~\cite[Theorem~2.5.7]{K80}).
On the other hand, a rigorous construction of
an optimally controlled process $\zeta$, such
as a $\beta$-$\gamma$ strategy, is rather
non-trivial.
In the context of this paper, we construct a
unique strong solution to the SDE (\ref{SDE}) when
the controlled process $\zeta$ is a $\beta$-$\gamma$
strategy (see Theorem~\ref{lem:beta-gamma} in
Section~\ref{sec:beta-gamma}).
Despite the central role that such strategies play
in stochastic impulse control, we are not aware of
any such rigorous SDE result.
Furthermore, this construction allows for a probabilistic
derivation of the optimal expected discounted
running reward as well as the optimal expected
discounted reward from control exdenditure
functionals (see (\ref{b-g-Rh}) and (\ref{b-g-adm-form})
in Theorem~\ref{lem:beta-gamma}).
The construction that we make can most easily be adapted
to derive the existence of strong solutions to optimally
controlled SDEs that arise in other stochastic
impulse control problems, even in dimensions higher
than one.

The paper is organised as follows.
Section~\ref{sec:pr-form} presents the precise formulation
of the control problem that we solve, including all of the
assumptions that we make.
In Section~\ref{sec:ODE}, we derive several results associated
with a linear ODE that we need for the solution to the
stochastic control problem we consider.
In Section~\ref{sec:beta-gamma}, we prove that
the SDE (\ref{SDE}) has a unique strong solution when
the controlled process $\zeta$ is a $\beta$-$\gamma$
strategy and we derive analytic expressions for certain
associated functionals using probabilistic techniques.
We derive the complete solution to the control
problem that we consider in Section~\ref{sec:solution}.
Finally, we present several examples illustrating the
assumptions that we make and the results that we
establish in Section~\ref{sec:ex}.

%===============================================
\section{Formulation of the stochastic control problem }
\label{sec:pr-form}

Fix a filtered probability space $\bigl( \Omega, \fcal, (\fcal_t),
\bbp \bigr)$ satisfying the usual conditions and carrying a
standard one-dimensional $(\fcal_t)$-Brownian motion $W$.
We consider a dynamical system, the uncontrolled stochastic
dynamics of which are modelled by the SDE
\ben
\di X_t = b(X_t) \, \di t + \sigma (X_t) \, \di W_t , \quad
X_0 = x  > 0 , \label{SDE-0}
\een
and we make the following assumption. 

\begin{ass} \label{A1} {\rm
The functions $b , \sigma :  [0 ,\infty[ \mbox{} \rightarrow
\bbr$ are locally Lipschitz continuous and $\sigma (x)
> 0$ for all $x>0$.
} \end{ass}
This assumption implies that  the scale function $p$
and the speed measure $m$ of the diffusion associated
with the SDE (\ref{SDE-0}), which are given by
\begin{gather}
p(1) = 0 \quad \text{and} \quad p' (x) = \exp \biggl( - 2
\int _1^x \frac{b(s)}{\sigma^2 (s)} \, \di s \biggr)
\label{p} \\
\text{and} \quad m(\di x) = \frac{2}{\sigma ^2(x) p' (x)}
\, \di x , \label{m}
\end{gather}
are well-defined.
Additionally, we make the following assumption on the
boundary classification of the diffusion associated with
(\ref{SDE-0}). 

\begin{ass} \label{A2} {\rm
The boundary point 0 is inaccessible while the boundary
point $\infty$ is natural.
} \end{ass}
The state space of the linear diffusion associated
with the SDE (\ref{SDE-0}) is the interval $\ical =
\mbox{} ]0, \infty[$.
Recall that the boundary point $p \in \{ 0, \infty \}$ of
$\ical$ is called {\em inaccessible\/} if
$\bbp_x (T_p < \infty) = 0$ for all $x \in \ical$
and {\em accessible\/} otherwise.
Furthermore, if the boundary $p$ is inaccessible,
then it is {\em natural\/} if
\be
\lim _{x \in \ical , \, x \rightarrow p} \bbp_x (T_y < t)
= 0 \quad \text{for all } y \in \ical \text{ and } t > 0
\ee
and {\em entrance\/} otherwise, namely, if
\be
\lim _{x \in \ical , \, x \rightarrow p} \bbp_x (T_y < t)
> 0 \quad \text{for some } y \in \ical \text{ and } t > 0
\ee
(e.g., see Revuz and Yor~\cite[Definition~VII.3.9]{RY}).
In these expressions, $T_y$ is the first hitting time
of the set $\{ y \}$, which is defined by
\ben
T_y = \inf \left\{ t \geq 0 \mid \ X_t = y \right\} ,
\quad \text{for } y > 0 . \label{T-hitting}
\een
In Borodin and Salminen~\cite[II.1.6]{BS}, an inaccessible
boundary point is called {\em not-exit\/}, while a
natural (resp., entrance) boundary point is
called {\em natural\/} (resp., {\em entrance-not-exit\/}).
Integral conditions for the classification of a boundary
point $p \in \{ 0, \infty \}$ of $\ical$ in terms of the
scale function $p$ and the speed measure $m$
can be found in this reference.

We next consider the stochastic control problem defined
by (\ref{SDE})--(\ref{Lambda}).

\begin{de} \label{def1} {\rm
The family of all admissible controlled strategies is the
set of all $(\fcal_t)$-adapted c\`{a}dl\`{a}g processes
$\zeta$ with increasing and piece-wise constant sample
paths such that the SDE (\ref{SDE}) has a unique
non-explosive  strong solution and
\ben
\bbe_x \Biggl[ \sum _{t \geq 0} e^{-\Lambda _t^\zeta}
{\bf 1} _{\{ \Delta \zeta_t > 0 \}} \Biggr] < \infty . \label{AC}
\een
} \end{de}

\begin{ass} \label{A3} {\rm
The discounting rate function $r$ is bounded and continuous.
Also, there exists $r_0>0$ such that $r(x) \geq r_0$ for all
$x \geq 0$.
} \end{ass} 

To complete the set of our assumptions, we consider
the operator $\lscr$ acting on $C^1$ functions with
absolutely continuous first-order derivatives that is defined by
\ben
\lscr w(x) = \frac{1}{2} \sigma^2 (x) w''(x) + b(x) w'(x) - r(x)
w(x) . \label{L}
\een
In the presence of Assumptions \ref{A1}, \ref{A2} and \ref{A3},
the second-order linear ODE $\lscr w(x) = 0$ has two
fundamental $C^2$ solutions $\varphi$ and $\psi$ such that
\begin{gather}
0 < \varphi (x) \quad \text{and} \quad \varphi' (x) < 0 \quad
\text{for all } x>0 , \label{phipsi} \\
0 < \psi (x) \quad \text{and} \quad \psi' (x) > 0 \quad
\text{for all } x>0 \\
\text{and} \quad
\lim _{x \downarrow 0} \varphi (x)  = \lim _{x \uparrow \infty}
\psi (x) = \infty . \label{limphipsi}
\end{gather}
If 0 is a natural boundary point, then
\ben
\lim _{x \downarrow 0} \frac{\varphi' (x)}{p' (x)} = -\infty
, \quad
\lim _{x \downarrow 0} \psi (x) = 0
\quad \text{and} \quad
\lim _{x \downarrow 0} \frac{\psi' (x)}{p' (x)} = 0 ,
\label{0-natural}
\een
while, if 0 is an entrance boundary point, then
\ben
\lim _{x \downarrow 0} \frac{\varphi' (x)}{p' (x)} > -\infty
, \quad
\lim _{x \downarrow 0} \psi (x) > 0
\quad \text{and} \quad
\lim _{x \downarrow 0} \frac{\psi' (x)}{p' (x)} = 0 .
\label{0-entrance}
\een
Symmetric results hold for the boundary point $\infty$
(e.g., see Borodin and Salminen~\cite[II.10]{BS}).

The functions $\varphi$ and $\psi$ admit the probabilistic
representations
\ben
\varphi (y) = \varphi (x) \bbe _y \bigl[ e^{-\Lambda _{T_x}}
\bigr] \quad \text{and} \quad
\psi (x) = \psi (y) \bbe _x \bigl[ e^{-\Lambda _{T_y}} \bigr]
\quad \text{for all } x<y
, \label{phi-psi}
\een
where $\Lambda$ is defined by (\ref{Lambda}) with $X$
in place of $X^\zeta$ and $T_y$ is defined by
(\ref{T-hitting}).

Furthermore, $\varphi$ and $\psi$ are such that
\ben
\varphi (x) \psi' (x) - \varphi' (x) \psi (x) = C p' (x) \quad
\text{for all } x > 0 , \label{C}
\een
where $C = \varphi (1) \psi' (1) - \varphi' (1) \psi (1)$
and $p$ is the scale function defined by (\ref{p}).
To simplify the notation, we also define
\ben
\Psi (x) = \frac{2 \psi (x)}{C \sigma^2 (x) p'(x)} \quad
\text{and} \quad \Phi (x) = \frac{2 \varphi (x)}
{C \sigma^2 (x) p'(x)} . \label{PhiPsi}
\een

Beyond involving standard integrability and growth
assumptions, the conditions in the following assumption
may appear involved.
However, they are standard in the relevant literature
and are satisfied by a wide range of problem data
choices (see Examples~\ref{ex1}-\ref{ex4} in
Section~\ref{sec:ex}).

\begin{ass} \label{A5}  {\rm
The following conditions hold true:
\smallskip

\noindent (i)
The function $h$ is continuous as well as bounded from below.
Also, the limit $\lim _{x \downarrow 0} h(x) / r(x)$ exists in
$\bbr$ and
\be
\bbe_x \biggl[ \int _0^\infty e^{-\Lambda_t} \bigl| h(X_t)
\bigr| \, \di t \biggr] < \infty .
\ee

\noindent (ii)
The function $k$ is absolutely continuous,
\begin{gather}
\int _0^1 \bigl| k(s) \bigr| \, \di s < \infty \quad
\text{and the function } x \mapsto \int _0^x k(s) \,
\di s \text{ is bounded from below} . \label{int-k-bfb}
\end{gather}
Furthermore,
\be
\bbe_x \biggl[ \int _0^\infty e^{-\Lambda_t} \, \Bigl| \lscr
\biggl( \int _0^{\boldsymbol{\cdot}} k(s) \, \di s \biggr)
(X_t) \Bigr| \, \di t \biggr] < \infty \quad \text{and} \quad
\limsup _{x \uparrow \infty} \frac{1}{\psi (x)} \int _0^x
k(s) \,  \di s \in \bbr_+ .
\ee

\noindent (iii)
If we define 
\ben
\Theta (x) = h(x) + \lscr \biggl( \int_0^{\boldsymbol{\cdot}}
k(s) \, \di s\biggr)(x) ,  \label{Theta}
\een
then $\Theta$ is continuous and there exists a constant
$\xi \in \mbox{} ]0, \infty[$ such that the restriction of
$\Theta / r$ in $]0, \xi[$ (resp., in $]\xi, \infty[$)
is strictly increasing (resp., strictly decreasing).
} \end{ass}

%%%%%%%%%%%%%%%%%%%%%%%%%%%%%%%%%%
\section{Results associated with a linear ODE}
\label{sec:ODE}

Unless stated otherwise, the results in this section
hold true if the coefficients of (\ref{SDE-0}) satisfy the
usual Engelbert and Schmidt conditions, rather than the
stronger Assumption~\ref{A1}, and the boundary
points $0$, $\infty$ are inaccessible.
We start by recalling some standard results
that we will need and can be found in, e.g., Lamberton
and Zervos~\cite[Section 4]{LZ}.
Consider a Borel measurable function $F: \mbox{}
]0, \infty[ \rightarrow \bbr$ such that
\ben
\bbe_x \biggl[ \int _0^\infty e^{-\Lambda_t} \bigl| F(X_t)
\bigr| \, \di t \biggr] < \infty  \quad \text{for all } x>0 , \label{IC-0}
\een
where $\Lambda$ is defined by (\ref{Lambda}) for
$X^\zeta = X$.
This integrability condition is equivalent to
\ben
\int _0^x \bigl| F(s) \bigr| \Psi (s) \, \di s + \int _x^\infty
\bigl| F(s) \bigr| \Phi (s) \, \di s  < \infty \quad \text{for all }
x>0 , \label{IC}
\een
where $\Phi$ and $\Psi$ are defined by (\ref{PhiPsi}).
Given such a function $F$, we define
\ben
R_F (x) = \bbe_x \biggl[ \int _0^\infty e^{-\Lambda_t}
F(X_t) \, \di t \biggr] , \quad \text{for } x>0 . \label{RF}
\een
The function $R_F$ admits the analytic presentation
\ben
R_F (x) = \varphi (x) \int _0^x F(s) \Psi (s) \, \di s +
\psi (x) \int _x^\infty F(s) \Phi (s) \, \di s \label{RF1}
\een
and satisfies the ODE $\lscr R_F + F = 0$.
Furthermore,
\ben
\lim _{x \downarrow 0} \frac{\bigl| R_F (x) \bigr|}
{\varphi (x)} = 0 \quad \text{and} \quad \lim
_{x \uparrow \infty} \frac{\bigl| R_F (x) \bigr|}{\psi (x)}
= 0 . \label{limitRF/phiRF/psi}
\een

Conversely, consider any function $f : \mbox{} ]0, \infty[
\rightarrow \bbr$ that is $C^1$ with absolutely continuous
first-order derivative and such that
\be
\bbe_x \biggl[ \int _0^\infty e^{-\Lambda_t} \bigl| \lscr
f(X_t) \bigr| \, \di t \biggr] < \infty , \quad \limsup
_{z \downarrow 0} \frac{\bigl| f(z) \bigr|}{\varphi (z)}
< \infty \quad \text{and} \quad \limsup _{z \uparrow \infty}
\frac{\bigl| f(z) \bigr|}{\psi (z)} < \infty .
\ee
Such a function is such that
\begin{gather}
\text{both of the limits }
\lim_{z \downarrow 0} \frac{f(z)}{\varphi (z)} \text{ and }
\lim _{z \uparrow \infty} \frac{f(z)}{\psi (z)} \text{ exist }
\label{f-lims} \\
\text{and} \quad
f(x) = \lim _{z \downarrow 0} \frac{f(z)}{\varphi (z)}
\varphi (x) - R_{\lscr f} (x) + \lim _{z \uparrow \infty}
\frac{f(z)}{\psi (z)} \psi (x) \quad \text{for all } x>0 .
\label{f-RLf}
\end{gather}

Part~(ii) of the following result will be important in
appreciating the role that the boundary classification
of 0 has on whether switching off the system might
be optimal (see Remark~\ref{rem:nat/entr} at the
end of Section~\ref{sec:solution}).
In general, (\ref{R(0)-nat}) is not true if 0 is an entrance
boundary point  (see (\ref{Rh(0)}) in Example~\ref{ex8}
in Section~\ref{sec:ex}).

\begin{lem} \label{lem:R_F}
Suppose that Assumptions~\ref{A1} and~\ref{A3} hold true.
Also, suppose that the boundary points 0 and $\infty$ of the
diffusion associated with the SDE (\ref{SDE-0}) are both
inaccessible.
Let $F$ be any Borel measurable function satisfying the
equivalent integrability conditions (\ref{IC-0}) and (\ref{IC}),
and consider the function $R_F$ defined by (\ref{RF})
and (\ref{RF1}).
The following statements hold true:
\smallskip

\noindent {\rm (i)}
Suppose that $F$ is bounded from below.
If $K$ is any constant such that $F(x) / r(x) \geq K$ for all
$x>0$, then $R_F (x) \geq K$ for all $x>0$.
\smallskip

\noindent {\rm (ii)}
If 0 is a natural boundary point, then
\ben
\liminf _{x \downarrow 0} \frac{F(x)}{r(x)} \leq
\liminf _{x \downarrow 0} R_F (x) \leq
\limsup _{x \downarrow 0} R_F (x) \leq
\limsup _{x \downarrow 0} \frac{F(x)}{r(x)} . \label{R(0)-nat}
\een
\end{lem}
{\bf Proof.}
%=======
Part~(i) of the lemma follows immediately from the
calculation
\be
\inf_{x>0} R_F (x) = \inf _{x>0} \bbe_x \biggl[ \int
_0^\infty e^{-\Lambda_t} F(X_t) \, \di t \biggr] \geq
\inf _{x>0} \frac{F(x)}{r(x)} \bbe_x \biggl[ \int _0^\infty
e^{-\Lambda_t} r(X_t) \, \di t \biggr] =
\inf _{x>0} \frac{F(x)}{r(x)} ,
\ee
where we have used the definition (\ref{Lambda})
of $\Lambda$.

To establish part~(ii) of the lemma suppose in what follows
that 0 is a natural boundary point.
Assuming that $\limsup _{x \downarrow 0} F(x) / r(x)
\in \bbr$, fix any $\varepsilon > 0$ and let
$x_\varepsilon > 0$ be any point such that
\be
\frac{F(x)}{r(x)} \leq \limsup _{x \downarrow 0} \frac{F(x)}{r(x)}
+ \varepsilon \quad \text{for all } x \in \mbox{} ]0,
x_\varepsilon] .
\ee
In view of (\ref{RF}), (\ref{RF1}), the definition (\ref{Lambda})
of $\Lambda$ and the second limit in (\ref{0-natural}),
we can see that
\begin{align}
\limsup _{x \downarrow 0} R_F (x) - \limsup _{x \downarrow 0}
& \frac{F(x)}{r(x)} - \varepsilon \nonumber \\
= \mbox{} & \limsup _{x \downarrow 0} \bbe_x \biggl[ \int
_0^\infty e^{-\Lambda_t} \biggl( \frac{F(X_t)}{r(X_t)} -
\limsup _{x \downarrow 0} \frac{F(x)}{r(x)} - \varepsilon \biggr)
r(X_t) \, \di t \biggr] \nonumber \\
= \mbox{} & \limsup _{x \downarrow 0} \Biggl( \varphi (x)
\int _0^x \biggl( \frac{F(s)}{r(s)} - \limsup _{x \downarrow 0}
\frac{F(x)}{r(x)} - \varepsilon \biggr) r(s) \Psi (s) \, \di s
\nonumber \\
& + \psi (x) \int _x^\infty \biggl( \frac{F(s)}{r(s)} -
\limsup _{x \downarrow 0} \frac{F(x)}{r(x)} - \varepsilon \biggr)
r(s) \Phi (s) \, \di s \Biggr) \nonumber \\
\leq \mbox{} & \lim _{x \downarrow 0} \psi (x) \int
_{x_\varepsilon}^\infty \biggl( \frac{F(s)}{r(s)} - \limsup
_{x \downarrow 0} \frac{F(x)}{r(x)} - \varepsilon \biggr)
r(s) \Phi (s) \, \di s \Biggr) \nonumber \\
= \mbox{} & 0 , \nonumber
\end{align}
which implies that $\limsup _{x \downarrow 0} R_F (x)
\leq \limsup _{x \downarrow 0} F(x) / r(x)$ because
$\varepsilon$ has been arbitrary.
Similarly, we can show that $\lim _{x \downarrow 0}
R_F (x) = - \infty$ if $\lim _{x \downarrow 0} F(x) / r(x)
= -\infty$, and the third inequality in (\ref{R(0)-nat})
follows.
Using similar arguments, we can establish the first
inequality in (\ref{R(0)-nat}).
\mbox{} \hfill $\Box$

\begin{lem} \label{lem-G}
Suppose that Assumptions~\ref{A1} and~\ref{A3} hold true,
suppose that the boundary points 0 and $\infty$ of the
diffusion associated with the SDE (\ref{SDE-0}) are both
inaccessible and consider any Borel measurable function
$F$ satisfying the equivalent integrability conditions
(\ref{IC-0}) and (\ref{IC}).
The function $G_F : \mbox{} ]0,\infty[ \mbox{}
\rightarrow \bbr$ defined by
\ben
G_F (x) := R_F (x) - \frac{R'_F (x)}{\psi' (x)} \psi (x) =
\frac{C p' (x)}{\psi' (x)} \int _0^x F(s) \Psi (s) \, \di s
\label{G}
\een
is such that
\ben
\liminf _{x \downarrow 0} \frac{F(x)}{r(x)} \leq
\liminf _{x \downarrow 0} G_F (x) \leq \limsup
_{x \downarrow 0} G_F (x) \leq \limsup _{x \downarrow 0}
\frac{F(x)}{r(x)} . \label{GFlim0}
\een
Furthermore, if the boundary point $\infty$ is natural, then
\ben
\liminf _{x \uparrow \infty}
\frac{F(x)}{r(x)} \leq \liminf _{x \uparrow \infty} G_F (x) 
\leq \limsup _{x \uparrow \infty} G_F (x) \leq \limsup
_{x \uparrow \infty} \frac{F(x)}{r(x)} . \label{GFlimoo}
\een
\end{lem}
{\bf Proof.}
%=======
We first note that the equality in (\ref{G}) follows immediately
from the definition (\ref{RF1}) of $R_F$ and the identity
(\ref{C}).
In view of (\ref{0-natural}) and (\ref{0-entrance}), the
assumption that the boundary point 0 is inaccessible
implies that
\ben
\lim _{x \downarrow 0} \frac{\psi' (x)}{p' (x)} = 0 .
\label{psi'/p'-lim0}
\een
This limit and the calculation
\be
\frac{\di}{\di x} \frac{\psi' (x)}{p' (x)} = \frac{2}{\sigma^2 (x) p'(x)}
\biggl( \frac{1}{2} \sigma^2 (x) \psi'' (x) + b(x) \psi' (x) \biggr)
= \frac{2  r(x) \psi (x)}{\sigma^2 (x) p'(x)} = C r(x) \Psi (x)
\ee
imply that
\ben
\int _0^x r(s) \Psi (s) \, \di s = \frac{\psi' (x)}{C p' (x)} .
\label{Psi-int}
\een
Similarly, the calculation
\ben
\frac{\di}{\di x} \frac{\varphi' (x)}{p' (x)} = C r(x) \Phi (x)
\label{(phi'/p')'}
\een
and the assumption that the boundary point $\infty$ is
inaccessible imply that
\ben
\int _x^\infty r(s) \Phi (s) \, \di s = - \frac{\varphi' (x)}{C p' (x)} .
\label{Phi-int}
\een
In view of (\ref{Psi-int}) and the expression of $G_F$ on the
right-hand side of (\ref{G}), we can see that
\begin{align}
G_F (x) & \geq \frac{C p' (x)}{\psi' (x)} \inf _{y<x} \frac{F(y)}{r(y)}
\int _0^x r(s) \Psi(s) \, \di s = \inf _{y<x} \frac{F(y)}{r(y)}
\nonumber \\
\text{and} \quad
G_F (x) & \leq \frac{C p' (x)}{\psi' (x)} \sup _{y<x} \frac{F(y)}{r(y)}
\int _0^x r(s) \Psi (s) \, \di s = \sup _{y<x} \frac{F(y)}{r(y)}
. \nonumber
\end{align}
These inequalities imply (\ref{GFlim0}).

Next, we additionally assume that $\infty$ is a natural
boundary point, which implies that $\lim _{x \uparrow \infty}
\psi' (x) / p' (x) = \infty$ (e.g., see Borodin and
Salminen~\cite[II.10]{BS}).
The expression of $G_F$ on the right-hand side of (\ref{G}),
the strict positivity of $\Psi$ and the identity (\ref{Psi-int})
imply that, given any $x>z>0$,
\begin{align}
\frac{C p' (x)}{\psi' (x)} \int _0^z F(s) & \Psi (s) \, \di s
+ \inf _{y>z} \frac{F(y)}{r(y)} \biggl( 1 - \frac{p' (x)}{\psi' (x)}
\frac{\psi' (z)}{p' (z)} \biggr) \nonumber \\
& = \frac{C p' (x)}{\psi' (x)} \int _0^z F(s) \Psi (s) \, \di s
+ \inf _{y>z} \frac{F(y)}{r(y)} \frac{C p' (x)}{\psi' (x)}
\int _z^x r(s) \Psi (s) \, \di s \nonumber \\
& \leq G_F (x)
\leq \frac{C p' (x)}{\psi' (x)} \int _0^z F(s) \Psi (s) \, \di s
+ \sup _{y>z} \frac{F(y)}{r(y)} \frac{C p' (x)}{\psi' (x)}
\int _z^x r(s) \Psi (s) \, \di s \nonumber \\
& = \frac{C p' (x)}{\psi' (x)} \int _0^z F(s)
\Psi (s) \, \di s + \sup _{y>z} \frac{F(y)}{r(y)} \biggl(
1 - \frac{p' (x)}{\psi' (x)} \frac{\psi' (z)}{p' (z)} \biggr) .
\nonumber
\end{align}
Combining these observations, we can see that
\be
\inf _{y>z} \frac{F(y)}{r(y)} \leq \liminf _{x \uparrow \infty}
G_F (x) \leq \limsup _{x \uparrow \infty} G_F (x) \leq 
\sup _{y>z} \frac{F(y)}{r(y)} \quad \text{for all } z >0 ,
\ee
and (\ref{GFlimoo}) follows.
\mbox{} \hfill $\Box$

\begin{lem} \label{R'/psi'=0?}
Suppose that Assumption~\ref{A1} and~\ref{A3} hold
true.
Also, suppose that the boundary points 0 and $\infty$
of the diffusion associated with the SDE (\ref{SDE-0})
are both inaccessible.
Given any Borel measurable function $F$ satisfying the
equivalent integrability conditions (\ref{IC-0}) and (\ref{IC}),
if the boundary point 0 (resp., $\infty$) is inaccessible,
then
\ben
\liminf _{x \downarrow 0} \frac{R'_F (x)}{\varphi' (x)}
\leq 0 \leq \limsup _{x \downarrow 0} \frac{R'_F (x)}{\varphi' (x)}
\quad \biggl( \text{resp., }
\liminf _{x \uparrow \infty} \frac{R'_F (x)}{\psi' (x)} \leq 0
\leq \limsup _{x \uparrow \infty} \frac{R'_F (x)}{\psi' (x)}
\biggr) . \label{R'_F/phi-psi'}
\een
Furthermore, if there exists $x_\dagger > 0$ {\rm (}resp.,
$x^\dagger > 0${\rm )} such that the restriction of $F/r$
in $]0, x_\dagger[$  {\rm (}resp., $]x^\dagger, \infty[${\rm )}
is a monotone function, then 
\ben
\lim_{x \downarrow 0} \frac{R'_F (x)}{\varphi' (x)} = 0 \quad
\biggl( \text{resp., } \lim_{x \uparrow \infty}
\frac{R'_F (x)}{\psi' (x)} = 0 \biggr) . \label{R'_F/phi-psi'-lims}
\een
\end{lem}
{\bf Proof.} 
%=======
To establish the very first inequality in (\ref{R'_F/phi-psi'}),
we argue by contradiction. 
To this end, we assume that
$\liminf _{x \downarrow 0} R'_F(x) / \varphi'(x) > 0$, which
implies that there exist $\varepsilon>0$ and $x_\varepsilon
>0$ such that
\be
\frac{R'_F (x)}{\varphi' (x)} > \varepsilon \ \Leftrightarrow
\  R'_F (x) < \varepsilon \varphi' (x) \quad \text{for all }
x \in \mbox{} ]0, x_\varepsilon[ .
\ee
However, this observation and the fact that  $\lim
_{x \downarrow 0} \varphi (x) = \infty$ imply that
$\lim _{x \downarrow 0} R_F (x) / \varphi (x) \geq \varepsilon$,
which contradicts (\ref{limitRF/phiRF/psi}).
The proof of the other inequalities in (\ref{R'_F/phi-psi'})
is similar.

To proceed further, we first note that (\ref{C}) and the fact that
$\lscr \varphi = \lscr \psi = 0$, where $\lscr$ is the differential
operator defined by (\ref{L}), imply that
\ben
\psi' (x) \varphi'' (x) - \varphi' (x) \psi'' (x) =
\frac{2C r(x)}{\sigma^2 (x)} p'(x) . \label{C-1}
\een
In view of this observation and the definition (\ref{RF1}) of
$R_F$, we can see that  the function $R'_F / \psi'$
is absolutely continuous with derivative
\begin{align}
\frac{\di}{\di x} \frac{R'_F (x)}{\psi' (x)} & =
\frac{2C r(x) p'(x)}{( \sigma (x) \psi' (x) )^2}
\biggl( \int _0^x F(s) \Psi (s) \, \di s - \frac{F(x)}{r(x)}
\frac{\psi' (x)}{Cp' (x)} \biggr) =: \frac{2C r(x) p'(x)}
{( \sigma (x) \psi' (x) )^2} Q_F (x) . \label{(R'/psi')'}
\end{align}
Now, suppose that there exists a point $x^\dagger > 0$
such that $F/r$ is monotone in $[x^\dagger, \infty[$.
Given any points $x_1 < x_2$ in $[x^\dagger, \infty[$,
we use (\ref{Psi-int}) to calculate
\begin{align}
Q_F (x_2) - Q_F (x_1) & = \int _{x_1}^{x_2} F(s) \Psi (s)
\, \di s - \frac{F(x_2)}{r(x_2)} \frac{\psi' (x_2)}{C p' (x_2)}
+ \frac{F(x_1)}{r(x_1)} \frac{\psi' (x_1)}{C p' (x_1)}
\nonumber \\
& = \int _{x_1}^{x_2} \biggl( \frac{F(s)}{r(s)}
- \frac{F(x_2)}{r(x_2)} \biggr) r(s) \Psi (s) \, \di s  +
\frac{\psi' (x_1)}{C p' (x_1)} \biggl( \frac{F(x_1)}{r(x_1)}
- \frac{F(x_2)}{r(x_2)} \biggr) \nonumber \\
& \begin{cases} \geq 0 , & \text{if } F/r \text{ is decreasing in }
[x^\dagger, \infty[ , \\ \leq 0 , & \text{if } F/r
\text{ is increasing in } [x^\dagger, \infty[ . \end{cases}
\label{Q}
\end{align}
Therefore, $Q_F$ is monotone in $[x^\dagger, \infty[$
and the limit $\lim _{x \uparrow \infty} Q_F (x)$ exists in
$[{-\infty}, \infty]$.
However, this observation and (\ref{(R'/psi')'}) imply that
there exists $\widetilde{x} \geq x^\dagger$ such that
$R'_F / \psi'$ is monotone in $[\widetilde{x}, \infty[$.
Therefore, the limit $\lim _{x \uparrow \infty} R'_F (x)
/ \psi' (x)$ exists, which, combined with the last two
inequalities in (\ref{R'_F/phi-psi'}), implies the
corresponding limit in (\ref{R'_F/phi-psi'-lims}).

Finally, we can establish the other limit in
(\ref{R'_F/phi-psi'-lims}) using symmetric arguments
and (\ref{(phi'/p')'}).
\mbox{} \hfill $\Box$

\bigskip

The following result will play a critical role in our
analysis.
Example~\ref{ex9} in Section~\ref{sec:ex} shows
that the point $\xlb$ introduced in part~(i) of the lemma
can be equal to $\infty$ if the sufficient conditions in
(\ref{xlb<oo-SC}) fail to be true.
Also, in contrast to the limit in (\ref{z_*}), Examples~\ref{ex5}
and~\ref{ex6} in Section~\ref{sec:ex} show that the limit
$\lim_{x \downarrow 0} R'_\Theta (x) / \psi' (x)$, which
characterises part~(iii) of the lemma, can take any value
in $]{-\infty}, \infty]$.

\begin{lem} \label{lem-z*}
Suppose that Assumption~\ref{A1} and~\ref{A3} hold true.
Also, suppose that the boundary points 0 and $\infty$ are
both inaccessible.
Given a function $\Theta$ satisfying the conditions
of Assumption~\ref{A5}.(iii), as well as  the equivalent
integrability conditions (\ref{IC-0}) and (\ref{IC}),
the following statements are true:
\smallskip

\noindent {\rm (i)}
There exists a unique $\xlb \in \mbox{} ]\xi, \infty]$
such that
\ben
\frac{\di}{\di x} \frac{R'_\Theta (x)}{\psi' (x)} \begin{cases} < 0
& \text{for all } x \in \mbox{} ]0, \xlb[ , \\ > 0 & \text{for all }
x \in \mbox{} ]\xlb, \infty[ , \end{cases}
\quad \text{and} \quad
\lim _{x \uparrow \infty} \frac{R'_\Theta (x)}{\psi' (x)} = 0 ,
\label{z_*}
\een
where we adopt the convention $]\infty, \infty[ \mbox{}
= \emptyset$.
\smallskip

\noindent {\rm (ii)}
$\xlb < \infty$ if and only if $\lim _{x \uparrow \infty}
Q_ \Theta (x) > 0$, where
\ben
Q_ \Theta (x) = \int _0^x \Theta (s) \Psi (s) \, \di s -
\frac{\Theta (x)}{r(x)} \frac{\psi' (x)}{Cp' (x)} .
\label{Q-Theta}
\een
In particular, this is the case if
\ben
\lim _{x \uparrow \infty} \frac{\Theta (x)}{r(x)} =
- \infty \quad \text{or} \quad
\lim _{x \downarrow 0} \frac{\Theta (x)}{r(x)}
> \lim _{x \uparrow \infty} \frac{\Theta (x)}{r(x)}
. \label{xlb<oo-SC}
\een

\noindent {\rm (iii)}
If $\xlb < \infty$ and we define
\ben 
\xub = \inf \biggl\{ s>0 \Bigm| \ \frac{R'_\Theta (s)}{\psi' (s)}
> \lim _{x \downarrow 0} \frac{R'_\Theta (x)}{\psi' (x)} \biggr\}
, \label{z^*}
\een
with the usual convention that $\inf \emptyset = \infty$,
then $\xub > \xlb$, 
\begin{gather}
\xub = \infty \quad \Leftrightarrow \quad \lim
 _{x \downarrow 0}  \frac{R'_\Theta (x)}{\psi' (x)} \geq 0
\label{z*-equiv} \\
\text{and} \quad
\lim _{x \downarrow 0} \frac{\Theta (x)}{r(x)} = -\infty \quad
\Rightarrow \quad \lim _{x \downarrow 0}
\frac{R'_\Theta (x)}{\psi' (x)} = \infty \quad \Rightarrow
\quad \xub = \infty . \label{z*-impl}
\end{gather}
\end{lem}
{\bf Proof.} 
%=======
The limit in (\ref{z_*}) follows from Lemma~\ref{R'/psi'=0?}
and the assumption that $\Theta / r$ is strictly decreasing
in $]\xi, \infty[$.
Using (\ref{Psi-int}) and (\ref{(R'/psi')'}), we can see that
\be
\frac{\di}{\di x} \frac{R'_\Theta (x)}{\psi' (x)} =
\frac{2C r(x) p'(x)}{( \sigma (x) \psi' (x) )^2}
\int _0^x \biggl( \frac{\Theta (s)}{r(s)} -
\frac{\Theta (x)}{r(x)} \biggr) r(s) \Psi (s) \, \di s =
\frac{2C r(x) p'(x)}{( \sigma (x) \psi' (x) )^2} Q_ \Theta
(x) . \label{(R'/psi')'Theta}
\ee
These expressions imply that
\be
\frac{\di}{\di x} \frac{R'_\Theta (x)}{\psi' (x)} < 0
\quad \text{and} \quad Q_ \Theta (x) < 0 \quad
\text{for all } x \leq \xi
\ee
because $\Theta / r$ is strictly increasing in $]0, \xi[$.
On the other hand, (\ref{Q}) for $F = \Theta$ implies
that $Q_\Theta$ is strictly increasing in $[\xi, \infty[$
because $\Theta / r$ is strictly decreasing in
$[\xi, \infty[$.
It follows that there exists a unique $\xlb \in \mbox{}
]\xi, \infty]$ such that the inequalities in (\ref{z_*})
hold true.
Furthermore, $\xlb < \infty$ if and only if
$\lim _{x \uparrow \infty} Q_\Theta (x) > 0$.

To establish the sufficient conditions in part~(ii) of the
lemma, we first use the integration by parts formula
and (\ref{Psi-int}) to observe that
\begin{align}
Q_ \Theta (x) & = \int _0^\xi \Theta (s) \Psi (s) \, \di s
- \frac{\Theta (\xi)}{r(\xi)} \frac{\psi' (\xi)}{Cp' (\xi)}
- \int _\xi^x \frac{\psi' (s)}{C p' (s)} \, \di
\frac{\Theta (s)}{r(s)}\nonumber \\
& \geq \int _0^\xi \Theta (s) \Psi (s) \, \di s -
\frac{\Theta (x)}{r(x)} \frac{\psi' (\xi)}{Cp' (\xi)}
\quad \text{for al } \xi < x . \label{Q-Theta-x>>0}
\end{align}
This inequality reveals that $\lim _{x \uparrow \infty}
Q_ \Theta (x) = \infty$ if $\lim _{x \uparrow \infty}
\Theta (x) / r(x) = - \infty$.

The identity (\ref{Psi-int}) implies that,
given any constant $K$,
\be
\int _0^x K r(s) \Psi (s) \, \di s - \frac{K r(x)}{r(x)}
\frac{\psi' (x)}{C p' (x)} = 0
\ee
Combining this observation with the definition of
$Q_\Theta$, we can see that $Q_\Theta =
Q_{\Theta + Kr}$.
If $\Theta / r$ satisfies the inequality in
(\ref{xlb<oo-SC}), then, for all $K$ such that
\be
- \lim _{x \downarrow 0} \frac{\Theta (x)}{r(x)} <
K < - \lim _{x \uparrow \infty} \frac{\Theta (x)}{r(x)} ,
\ee
there exists $\eta (K) \in \mbox{} ]\xi, \infty[$
such that
\be
\Theta \bigl( \eta (K) \bigr) + K r \bigl( \eta (K) \bigr) =0
\quad \text{and} \quad
Q_ {\Theta + Kr} \bigl( \eta (K) \bigr) = \int _0^{\eta (K)}
\bigl( \Theta (s) + K r(s) \bigr) \Psi (s) \, \di s > 0 .
\ee
It follows that
\be
\lim _{x \uparrow \infty} Q_\Theta (x) =
\lim _{x \uparrow \infty} Q_{\Theta + Kr} (x) > 0 ,
\ee
thanks to the fact that $Q_\Theta$ is strictly increasing
in $[\xi, \infty[$.

The equivalence (\ref{z*-equiv}) follows immediately
from (\ref{z_*}) and the definition (\ref{z^*}) of $\xub$.
To establish the implications in (\ref{z*-impl}), we first
note that (\ref{(phi'/p')'}) implies that the function
$\varphi' / p'$ is strictly increasing, so the limit
$\lim _{x \downarrow 0} \varphi' (x) / p'(x)$ exists in
$[-\infty, 0[$.
Therefore,
\ben
\lim _{x \downarrow 0} \frac{\psi' (x)}{\varphi' (x)} =
\lim _{x \downarrow 0} \frac{\psi' (x)}{p' (x)}
\lim _{x \downarrow 0} \frac{p' (x)}{\varphi' (x)} = 0 ,
\label{psi'/phi'-lim0}
\een
where we have also used (\ref{psi'/p'-lim0}).
Using the first of these two observations, the definition
(\ref{RF1}) of $R_\Theta$, (\ref{Phi-int}), (\ref{C-1}) and
integration by parts, we can see that, if
$\lim _{x \downarrow 0} \Theta (x) / r(x) = - \infty$,
then
\begin{align}
\lim _{x \downarrow 0} \frac{( \sigma (x) \varphi' (x) )^2}
{2C r(x) p'(x)} \frac{\di}{\di x} \frac{R'_\Theta (x)}
{\varphi' (x)} & = - \lim _{x \downarrow 0} \biggl(
\int _x^\infty \Theta (s) \Phi (s) \, \di s + \frac{\Theta (x)}{r(x)}
\frac{\varphi' (x)}{C p' (x)} \biggr) \nonumber \\
& = - \int _1^\infty \Theta (s) \Phi (s) \, \di s - \frac{\Theta (1)}{r(1)}
\frac{\varphi' (1)}{C p' (1)} + \lim _{x \downarrow 0}
\int_x^1 \frac{\varphi' (s)}{C p' (s)} \, \di \frac{\Theta (s)}{r(s)}
\nonumber \\
& = - \infty . \nonumber
\end{align}
On the other hand, we use (\ref{C-1}) to calculate
\be
\frac{\di}{\di x} \frac{\psi' (x)}{\varphi' (x)} = -
\frac{2C r(x) p'(x)}{( \sigma (x) \varphi' (x) )^2} .
\ee
In view of (\ref{R'_F/phi-psi'-lims}) and
(\ref{psi'/phi'-lim0}), these calculations and
L'\,H\^{o}pital's formula imply that
\be
\lim _{x \downarrow 0} \frac{R'_\Theta (x)}{\psi' (x)}
= \lim _{x \downarrow 0} \frac{\frac{\di}{\di x}
\frac{R'_\Theta  (x)}{\varphi' (x)}}{\frac{\di}{\di x}
\frac{\psi' (x)}{\varphi' (x)}} = \infty .
\ee
The implications in (\ref{z*-impl}) follow from
this analysis and the definition (\ref{z^*}) of $\xub$.
\mbox{} \hfill $\Box$
\bigskip

%%%%%%%%%%%%%%%%%%%%%%%%%%%%%%%%%%
\section{The ``$\pmb \beta$-$\pmb \gamma$'' strategy}
\label{sec:beta-gamma}

In this section, we consider the $\beta$-$\gamma$
strategy that is characterised by two points
$0 < \gamma < \beta < \infty $ and takes the following
form.
If the state process  takes any value $x \geq \beta$,
the controller pushes it in an impulsive way down to
the level $\gamma$.
For as long as the state process takes values inside the
interval $]0, \beta[$, the controller waits and takes
no action.
Accordingly, such a strategy is characterised by a
controlled process $\zeta$ such that
\ben
\Delta \zeta_t = \bigl( X_{t-}^\zeta - \gamma \bigr) {\bf 1}
_{\{ X_{t-}^\zeta \geq \beta \} } \quad \text{for all } t \geq 0 ,
\label{Xzeta-BG}
\een
where $X^\zeta$ is the associated solution to the SDE
(\ref{SDE}).

%==============================================
\begin{thm} \label{lem:beta-gamma}
Suppose that Assumptions~\ref{A1} and \ref{A3} hold
true.
Also, suppose that the boundary points 0 and $\infty$
of the diffusion associated with the uncontrolled SDE
(\ref{SDE-0}) are both inaccessible.
Given any points $\gamma < \beta$ in $]0, \infty[$, there
exists a controlled process $\zeta = \zeta (\beta, \gamma)$
that is admissible in the sense of Definition~\ref{def1}
and is such that (\ref{Xzeta-BG}) holds true.
Furthermore, given any $x \in \mbox{} ]0, \beta[$,
\ben
\bbe_x \Biggl[ \int _0^\infty e^{-\Lambda_t^\zeta}
h (X_t^\zeta) \, \di t \biggr] = R_h (x) +
\frac{\psi (x)}{\psi (\beta) - \psi (\gamma)}
\bigl( R_h (\gamma) - R_h (\beta) \bigr) \label{b-g-Rh}
\een
and
\ben
\bbe_x \Biggl[ \sum _{t \geq 0} e^{-\Lambda _t^\zeta}
{\bf 1} _{\{ \Delta \zeta_t > 0 \}} \Biggr] =
\frac{\psi (x)}{\psi (\beta) - \psi (\gamma)} .
\label{b-g-adm-form}
\een
\end{thm}
{\bf Proof.}
%=======
We start with a recursive construction of the required
process $\zeta$ and its associated solution to the
SDE (\ref{SDE-0}).
To this end, we first consider any initial state $x \in
\mbox{} ]0, \beta[$, we denote by $X^1$ the solution
to the uncontrolled SDE (\ref{SDE-0}) and we define
\ben
\tau_1 = \inf \bigl\{ t \geq 0 \mid \ X_t^1 \geq \beta \bigr\}
\quad \text{and} \quad
\zeta _t^1 = (\beta - \gamma) {\bf 1} _{\{ \tau_1 \leq t \}} .
\label{zeta1}
\een
Given $\ell \geq 1$, suppose that we have determined
$X^j$, $\tau_j$ and $\zeta^j$, for $j = 1, \ldots, \ell$.

The process $\widetilde{W}^{\ell+1}$
defined by $\widetilde{W}_t^{\ell+1} = \bigl( W_{\tau_\ell + t}
- W_{\tau_\ell} \bigr) {\bf 1} _{\{ \tau_\ell < \infty \}}$ is a standard
$(\fcal _{\tau_\ell +t})$-Brownian motion  that is independent
of $\fcal _{\tau_\ell}$ under the conditional
probability measure $\bbp ( \cdot \mid \tau_\ell < \infty)$
(see Revuz and Yor~\cite[Exercise IV.3.21]{RY}).
We denote by $\widetilde{X}^{\ell+1}$ the unique
solution to the uncontrolled SDE (\ref{SDE-0}) with
$\widetilde{X} _0^{\ell+1} = \gamma$ that is driven by the
Brownian motion $\widetilde{W}^{\ell+1}$ and is defined
on the probability space
$\bigl( \Omega, \fcal, (\fcal _{\tau_\ell +t}), \bbp (\cdot
\mid \tau_\ell < \infty) \bigr)$.
Since $(t - \tau_\ell)^+$ is an $(\fcal _{\tau_\ell +t})$-stopping
time for all $t \geq 0$,
\be
\tau_\ell + (t - \tau_\ell)^+ = t \vee \tau_\ell
\quad \text{and} \quad
\widetilde{W} _{(t - \tau_\ell)^+}^{\ell+1} = \bigl(
W_{t \vee \tau_\ell} - W_{\tau_\ell} \bigr) {\bf 1}
_{\{ \tau_\ell < \infty \}} ,
\ee
we can see that, on the event $\{ \tau_\ell < \infty \}$,
\begin{align}
\widetilde{X} _{(t - \tau_\ell)^+}^{\ell+1} & = \gamma + \int
_0^{(t - \tau_\ell)^+} b \bigl( \widetilde{X} _s^{\ell+1} \bigr)
\, \di s + \int _0^{(t - \tau_\ell)^+} \sigma \bigl( \widetilde{X}
_s^{\ell+1} \bigr) \, \di \widetilde{W} _s^{\ell+1}
\nonumber \\
& = \gamma + \int _0^t b \Bigl( \widetilde{X}
_{(s - \tau_\ell)^+}^{\ell+1} \Bigr) \, \di (s - \tau_\ell)^+
+ \int _0^t \sigma \Bigl( \widetilde{X}
_{(s - \tau_\ell)^+}^{\ell+1} \Bigr) \, \di \widetilde{W}
_{(s - \tau_\ell)^+}^{\ell+1} \nonumber \\
& = \gamma + \int _{\tau_\ell}^{t \vee \tau_\ell}
b \Bigl( \widetilde{X} _{(s - \tau_\ell)^+}^{\ell+1} \Bigr)
\, \di s + \int _{\tau_\ell}^{t \vee \tau_\ell} \sigma \Bigl(
\widetilde{X} _{(s - \tau_\ell)^+}^{\ell+1} \Bigr) \, \di W_s ,
\nonumber
\end{align}
where we have time changed the Lebesgue as well as
the It\^{o} integral (see Revuz and
Yor~\cite[Propositions~V.1.4, V.1.5]{RY}).
It follows that, if we define
\ben
\overline{X}_t^{\ell+1} = \widetilde{X} _{(t - \tau_\ell)^+}^{\ell+1}
{\bf 1} _{\{ \tau_\ell < \infty \}} , \quad \text{for } t \geq 0 ,
\label{Xbar-ell}
\een
then 
\ben
\overline{X} _t^{\ell+1} = \gamma + \int
_{\tau_\ell}^{t \vee \tau_\ell} b \bigl( \overline{X} _s^{\ell+1}
\bigr) \, \di s + \int _{\tau_\ell}^{t \vee \tau_\ell} \sigma \bigl(
\overline{X} _s^{\ell+1} \bigr) \, \di W_s . \label{X-ell+1}
\een
Furthermore, we define
\begin{gather}
X _t^{\ell+1} = X_t^\ell {\bf 1} _{\{ t < \tau_\ell \}} +
\overline{X} _t^{\ell+1} {\bf 1} _{\{ \tau_\ell \leq t \}}
, \\
\tau_{\ell+1} = \inf \bigl\{ t > \tau_\ell \mid \ X_t^{\ell+1} \geq
\beta \bigr\} \quad \text{and} \quad
\zeta _t^{\ell+1} = \zeta _t^\ell + (\beta - \gamma) {\bf 1}
_{\{ \tau_{\ell+1} \leq t \}} . \label{zeta-ell+1}
\end{gather}
Also, we note that
\ben
\tau_{\ell+1} - \tau_\ell = \widetilde{T} _\beta^{\ell+1}
:= \inf \bigl\{ t \geq 0 \mid \ \widetilde{X} _t^{\ell+1} \geq
\beta \bigr\} .
\een

Given the recursive construction we have just considered,
we define
\ben
X_t^\zeta = \sum _{\ell = 0}^\infty X _t^{\ell+1} {\bf 1}
_{\{ \tau_\ell \leq t < \tau_{\ell+1} \}}
\quad \text{and} \quad
\zeta_t = \sum _{\ell = 0}^\infty \zeta _t^{\ell+1} {\bf 1}
_{\{ \tau_\ell \leq t < \tau_{\ell+1} \}} . \label{SDE-sol}
\een
In view of (\ref{X-ell+1})--(\ref{zeta-ell+1}), the process
$X^\zeta$ given by (\ref{SDE-sol}) provides the unique
solution to the SDE (\ref{SDE}) for $\zeta$ being as in
(\ref{SDE-sol}).
Furthermore, these processes are such that
(\ref{Xzeta-BG}) holds true.
In the case that arises when the initial state $x \geq \beta$,
the only modification of the arguments above involves
$X^1$ being the solution to the uncontrolled SDE
(\ref{SDE-0}) for $x = \gamma$ and $\zeta^1$ being
the same as in (\ref{zeta1}) translated by adding
the constant $x - \gamma$ to it.

We next establish (\ref{b-g-adm-form}), which implies
the admissibility condition (\ref{AC}).
The process $\widetilde{X}^{\ell+1}$ introduced at the
beginning of the proof is independent of $\fcal _{\tau_\ell}$
under the conditional probability measure $\bbp (\cdot
\mid \tau_\ell < \infty)$ and its distribution under $\bbp
(\cdot \mid \tau_\ell < \infty)$ is the same as the distribution
of the solution $X$ to the uncontrolled SDE (\ref{SDE-0})
with initial state $X_0 = \gamma$ under $\bbp$.
In particular,
\be
\bbe ^{\bbp (\cdot \mid \tau_\ell < \infty)} \Bigl[ F \bigl(
\widetilde{X}^{\ell+1} \bigr) \Bigr] = \bbe _\gamma \bigl[
F(X) \bigr]
\ee
for every bounded measurable functional $F$ mapping
continuous functions on $\bbr_+$ to $\bbr_+$, where
we denote by $\bbe ^{\bbp (\cdot \mid \tau_\ell < \infty)}$
expectations computed under the conditional probability
measure $\bbp (\cdot \mid \tau_\ell < \infty)$.
In view of these observations and the definition of conditional
expectation,
\ben
\bbe _\gamma \Bigl[ F \bigl( \widetilde{X}^{\ell+1} \bigr)
\, \big| \, \fcal _{\tau_\ell} \Bigr] {\bf 1} _{\{ \tau_\ell < \infty \}}
= \bbe _\gamma \bigl[ F (X) \bigr] {\bf 1} _{\{ \tau_\ell < \infty \}}
. \label{CondExp}
\een
To see this claim, we first note that
the Radon-Nikodym derivative of $\bbp (\cdot \mid
\tau_\ell < \infty)$ with respect to $\bbp$ is given by
\be
\frac{\di \bbp (\cdot \mid \tau_\ell < \infty)}{\di \bbp}
= \frac{1}{\bbp (\tau_\ell < \infty)} {\bf 1}
_{\{ \tau_\ell < \infty \}} .
\ee
Given any event $\Gamma \in \fcal _{\tau_\ell}$,
\begin{align}
\frac{1}{\bbp (\tau_\ell < \infty)} \bbe _\gamma \Bigl[
\bbe _\gamma \bigl[ F(X) \bigr] {\bf 1} _{\{ \tau_\ell < \infty \}}
{\bf 1} _\Gamma \Bigr] & = \bbe _\gamma \bigl[ F(X)
\bigr] \frac{1}{\bbp (\tau_\ell < \infty)} \bbe _\gamma
\bigl[ {\bf 1} _{\{ \tau_\ell < \infty \} \cap \Gamma} \bigr]
\nonumber \\
& = \bbe ^{\bbp (\cdot \mid \tau_\ell < \infty)} \Bigl[ F
\bigl( \widetilde{X}^{\ell+1} \bigr) \Bigr] \bbe
^{\bbp (\cdot \mid \tau_\ell < \infty)} \bigl[ {\bf 1}
_\Gamma \bigr] \nonumber \\
& = \bbe ^{\bbp (\cdot \mid \tau_\ell < \infty)} \Bigl[ F
\bigl( \widetilde{X}^{\ell+1} \bigr) {\bf 1} _\Gamma \Bigr]
\nonumber \\
& = \frac{1}{\bbp (\tau_\ell < \infty)} \bbe _\gamma \Bigl[
F \bigl( \widetilde{X}^{\ell+1} \bigr) {\bf 1}
_{\{ \tau_\ell < \infty \}} {\bf 1} _\Gamma \Bigr] , \nonumber
\end{align}
and (\ref{CondExp}) follows.

In view of (\ref{Xbar-ell})--(\ref{CondExp}), we can
see that
\begin{align}
\bbe _x \bigl[ e^{-\Lambda _{\tau_{\ell+1}}^\zeta} \bigr] & =
\bbe _x \biggl[ e^{-\Lambda _{\tau_\ell}^\zeta} \bbe _\gamma
\biggl[ \exp \bigg( - \int _{\tau_\ell}^{\tau_{\ell+1}} r(X^\zeta_u)
\, \di u \biggr) \, \Big| \, \fcal _{\tau_\ell} \biggr] {\bf 1}
_{\{ \tau_\ell < \infty \}} \biggr] \nonumber \\
& = \bbe _x \biggl[ e^{-\Lambda _{\tau_\ell}^\zeta} \bbe
_\gamma \biggl[ \exp \bigg( - \int _{\tau_\ell}^{\tau_{\ell+1}}
r \bigl( \overline{X} _u^{\ell+1} \bigr) \, \di u \biggr)
\, \Big| \, \fcal _{\tau_\ell} \biggr] {\bf 1} _{\{ \tau_\ell < \infty \}}
\biggr] \nonumber \\
& = \bbe _x \biggl[ e^{-\Lambda _{\tau_\ell}^\zeta} \bbe
_\gamma \biggl[ \exp \bigg( - \int _0^{\widetilde{T} _\beta^{\ell+1}}
r \bigl( \widetilde{X}_u^{\ell+1} \bigr) \, \di u \biggr) \, \Big|
\, \fcal _{\tau_\ell} \biggr] {\bf 1} _{\{ \tau_\ell < \infty \}}
\biggr] \nonumber \\
& = \bbe _x \biggl[ e^{-\Lambda _{T_\ell}^\zeta} \bbe
_\gamma \biggl[ \exp \bigg( - \int _0^{T_\beta} r(X_u)
\, \di u \biggr) \biggr] {\bf 1} _{\{ \tau_\ell < \infty \}} \biggr]
\nonumber \\
& = \bbe _x \bigl[ e^{-\Lambda _{\tau_\ell}^\zeta} \bigr]
\bbe _\gamma \bigl[ e^{-\Lambda _{T_\beta}} \bigr] ,
\nonumber
\end{align}
where $\Lambda$ is defined by (\ref{Lambda}) with $X$
in place of $X^\zeta$ and $T_\beta$ is defined as in
(\ref{T-hitting}).
Given any $x \in \mbox{} ]0, \beta[$, we iterate this result
and use (\ref{phi-psi}) to obtain
\ben
\bbe _x \bigl[ e^{-\Lambda _{\tau_{\ell+1}}^\zeta} \bigr]
= \bbe _x \bigl[ e^{-\Lambda _{T_\beta}} \bigr]
\Bigl( \bbe _\gamma \bigl[ e^{-\Lambda _{T_\beta}}
\bigr] \Bigr) ^\ell = \frac{\psi (x)}{\psi (\beta)} \biggl(
\frac{\psi (\gamma)}{\psi (\beta)} \biggr) ^\ell .
\label{E[eLam]}
\een
It follows that
\begin{align}
\bbe_x \Biggl[ \sum _{t \geq 0} e^{-\Lambda _t^\zeta}
{\bf 1} _{\{ \Delta \zeta_t > 0 \}} \Biggr] & = \bbe_x \Biggl[
\sum _{\ell = 1}^\infty e^{-\Lambda _{\tau_\ell}^\zeta}
\Biggr] \nonumber \\
& = \sum _{\ell = 1}^\infty \bbe_x \bigl[
e^{-\Lambda _{\tau_\ell}^\zeta} \bigr]
= \frac{\psi (x)}{\psi (\beta)} \sum _{\ell = 0}^\infty
\biggl( \frac{\psi (\gamma)}{\psi (\beta)} \biggr) ^\ell
= \frac{\psi (x)}{\psi (\beta) - \psi (\gamma)} ,
\nonumber
\end{align}
which establishes (\ref{b-g-adm-form}).

To show (\ref{b-g-Rh}), we consider any $x \in \mbox{}
]0, \beta[$ and we use (\ref{Xbar-ell})--(\ref{CondExp})
as well as (\ref{E[eLam]}) to derive the expression
\begin{align}
\bbe_x \biggl[ & \int _{\tau_\ell}^{\tau_{\ell+1}}
e^{-\Lambda_t^\zeta} h (X_t^\zeta) \, \di t \biggr]
\nonumber \\
& = \bbe_x \biggl[ e^{-\Lambda _{\tau_\ell}^\zeta}
\bbe _\gamma \biggl[ \int _{\tau_\ell}^{\tau_{\ell+1}}
\exp \bigg( - \int _{\tau_\ell}^t r \bigl( \overline{X}
_u^{\ell+1} \bigr) \, \di u \biggr) h (\overline{X} _t^{\ell+1})
\, \di t \, \Big| \, \fcal _{\tau_\ell} \biggr] \biggr]
\nonumber \\
& = \bbe _x \bigl[ e^{-\Lambda _{\tau_\ell}^\zeta} \bigr]
\bbe _\gamma \biggl[ \int _0^{T_\beta} e^{-\Lambda_t}
h(X_t) \, \di t \biggr]
= \bbe _x \bigl[ e^{-\Lambda _{\tau_\ell}^\zeta} \bigr]
\Bigl( R_h (\gamma) - \bbe _\gamma \bigl[
e^{-\Lambda _{T_\beta}} \bigr] R_h (\beta) \Bigr)
\nonumber \\
& =  \frac{\psi (x)}{\psi (\beta)} \biggl(
\frac{\psi (\gamma)}{\psi (\beta)} \biggr) ^{\ell-1}
\biggl( R_h (\gamma) - \frac{\psi (\gamma)}{\psi (\beta)}
R_h (\beta) \biggr) . \nonumber
\end{align}
Similarly, we can show that
\be
\bbe_x \biggl[ \int _0^{\tau_1} e^{-\Lambda_t^\zeta}
h (X_t^\zeta) \, \di t \biggr] = R_h (x) -
\frac{\psi (x)}{\psi (\beta)} R_h (\beta) .
\ee
Recalling the assumption that $h$ is bounded from
below, we can use the monotone convergence theorem
and these results to obtain
\begin{align}
\bbe_x \biggl[ & \int _0^\infty e^{-\Lambda_t^\zeta} h
(X_t^\zeta) \, \di t \biggr] \nonumber \\
& = \bbe_x \biggl[ \int _0^{\tau_1} e^{-\Lambda_t^\zeta}
h (X_t^\zeta) \, \di t \biggr] + \sum _{\ell=1}^\infty
\bbe_x \biggl[ \int _{\tau_\ell}^{\tau_{\ell+1}}
e^{-\Lambda_t^\zeta} h (X_t^\zeta) \, \di t \biggr]
\nonumber \\
& = R_h (x) - \frac{\psi (x)}{\psi (\beta)} R_h (\beta)
+ \frac{\psi (x)}{\psi (\beta)}
\biggl( R_h (\gamma) - \frac{\psi (\gamma)}{\psi (\beta)}
R_h (\beta) \biggr) \sum _{\ell=1}^\infty
\biggl( \frac{\psi (\gamma)}{\psi (\beta)} \biggr) ^{\ell-1}
, \nonumber
\end{align}
which proves (\ref{b-g-Rh}).
\mbox{} \hfill $\Box$

%%%%%%%%%%%%%%%%%%%%%%%%%%%%%%%%%%
\section{The solution to the control problem}
\label{sec:solution}

We will solve the control problem we have considered by
deriving a $C^1$ with absolutely continuous first-order
derivative function $w: \mbox{} ]0, \infty[\mbox{}
\rightarrow \bbr$ that satisfies the HJB equation
\ben
\max \Biggl\{ \lscr w(x) + h(x) , \ -c + \sup _{z \in [0,x[}
\int _{x-z}^x \bigl( k(s) - w'(s) \bigr) \, \di s \Biggr\} = 0
, \label{HJB}
\een
Lebesgue-a.e.\ in $]0,\infty[$.
Given such a solution, the optimal strategy can be characterised
as follows.
The controller should wait and take no action for as long
as the state process $X$ takes values in the interior of the
set in which the ODE
\be
\lscr w(x) + h(x) = 0
\ee
is satisfied and should take immediate action with an
impulse in the negative direction if the state process
takes values in the set of all points $x>0$ such that
\be
- c + \sup _{z \in [0,x[} \int _{x-z}^x \bigr( k(s) - w'(s)
\bigr) \, \di s  = 0 .
\ee

We first consider the possibility for a $\beta$-$\gamma$
strategy with $\gamma < \beta$ in $]0, \infty[$ to
be optimal.
The optimality of such a strategy is associated with
a solution $w$ to the HJB equation (\ref{HJB})
such that
\begin{gather}
\lscr w(x) + h(x) = 0 , \quad \text{for } x \in \mbox{}
]0, \beta[ , \label{OS1} \\
\text{and} \quad w(x) = w(\gamma) + \int _\gamma^x
k(s) \, \di s - c , \quad \text{for } x \in [\beta, \infty[ .
\label{OS2}
\end{gather}
To determine such a solution $w$, we first consider
the so-called ``principle of smooth fit", which requires
that $w'$ should be continuous, in particular, at the
free-boundary point $\beta$.
This condition suggests the free-boundary equation
\ben
\lim _{x \uparrow \beta} w'(x) = k(\beta) . \label{condb}
\een
Next we consider the inequality 
\be
- c + \sup _{z \in [0,x[} \int _{x-z}^x \bigr( k(s) - w'(s)
\bigr) \, \di s \leq 0 , 
\ee
which is associated with impulsive action.
For $x = \beta$ and $z = \beta - u$, we can see that
this implies that
\be
- c + \int _u^\beta \bigr( k(s) - w'(s) \bigr) \, \di s
\leq 0 \quad \text{for all } u \in \mbox{} ]0, \beta ] .
\ee
This inequality and the identity 
\ben
- c + \int _\gamma^\beta \bigr( k(s) - w'(s) \bigr)
\, \di s = 0,  \label{wgabe}
\een
which follows from (\ref{OS2}), can both be true if
and only if the function $u \mapsto \int _u^\beta
\bigr( k(s) - w'(s) \bigr) \, \di s$ has a local maximum
at $\gamma$.
This observation gives rise to the free-boundary
condition
\ben
w' (\gamma) = k(\gamma) . \label{condg}
\een

Every solution to (\ref{OS1}) that can satisfy the
so-called ``transversality condition'', which is
required for a solution $w$ to the HJB equation
to identify with the control problem's value function,
is given by
\ben
w(x) = R_h (x) + A \psi (x) , \label{w-cont}
\een
for some constant $A$, where $R_h$ is given by
(\ref{RF}) and (\ref{RF1}) for $F=h$.
In view of the definition (\ref{Theta}) of $\Theta$ in
Assumption~\ref{A5}, the expression of $R_\Theta$
as in (\ref{RF1}) and the representation (\ref{f-RLf}),
we can see that
\ben
R_h (x) = R_\Theta (x) + \int _0^x k(s) \, \di s - K_\infty
\psi (x) , \label{Rh-RTheta}
\een
where
\ben
K_\infty = \lim _{x \uparrow \infty} \frac{1}{\psi (x)}
\int _0^x k(s) \, \di s \in \bbr _+ . \label{Koo}
\een
Note that the limit $K_\infty$ indeed exists in $\bbr_+$,
thanks to the last condition in Assumption~\ref{A5}.(ii)
and (\ref{f-lims}).
The identity (\ref{Rh-RTheta}) implies that (\ref{w-cont})
is equivalent to
\be
w(x) = R_\Theta (x) + \int _0^x k(s) \, \di s + (A
- K_\infty) \psi (x) .
\ee
Therefore, the solution to (\ref{OS1}) that
satisfies the boundary condition (\ref{condb}) is given  by
\begin{align}
w(x) & = \int _0^x k(s) \, \di s + R_\Theta (x) -
\frac{R_\Theta' (\beta)}{\psi' (\beta)} \psi (x) \nonumber \\
& = R_h (x) + \biggl( K_\infty - \frac{R_\Theta' (\beta)}
{\psi' (\beta)} \biggr) \psi (x), \quad \text{for } x \in
\mbox{} ]0, \beta[ . \label{OS3}
\end{align}
Furthermore, the boundary conditions (\ref{condg})
and (\ref{wgabe}) are equivalent to
\ben
\frac{R_\Theta' (\gamma)}{\psi' (\gamma)} =
\frac{R_\Theta' (\beta)}{\psi' (\beta)} \quad \text{and}
\quad F(\gamma, \beta) = - c , \label{BC}
\een
respectively, where
\ben
F(\gamma, \beta) := G_\Theta (\beta) - G_\Theta (\gamma)
= \int _\gamma^\beta \biggl( \frac{R_\Theta' (s)}{\psi' (s)} -
\frac{R_\Theta' (\beta)}{\psi' (\beta)} \biggr) \psi' (s) \, \di s ,
\label{F}
\een
and $G_\Theta$ is defined by (\ref{G}) in
Lemma~\ref{lem-G}.

The following result is about the solvability of the system
of equations given by (\ref{BC}) for the unknowns
$\gamma$ and $\beta$.
Note that Lemma~\ref{lem-z*}.(i) implies that a pair
$0 \leq \gamma < \beta < \infty$ satisfying the first equation
in (\ref{BC}) might exist only if $\xlb < \infty$.

\begin{lem} \label{lem:FB}
Consider the stochastic control problem formulated in
Section~\ref{sec:pr-form} and suppose that the
point $\xlb$ introduced in Lemma~\ref{lem-z*}.(i) is
finite.
There exist a unique strictly decreasing function
$\gamma^\star: \mbox{} ]0, c^\star[ \mbox{}
\rightarrow \mbox{} ]0, \xlb[$ and a unique strictly
increasing function $\beta^\star: \mbox{} ]0, c^\star[
\mbox{} \rightarrow \mbox{} ]\xlb, \xub[$, where
$c^\star > 0$ is defined by (\ref{c*}) in the proof below
and $\xlb$, $\xub$ are as in Lemma~\ref{lem-z*},
such that
\ben
\frac{R_\Theta' (x)}{\psi' (x)} -
\frac{R_\Theta' \bigl( \beta^\star (c) \bigr)}
{\psi' \bigl( \beta^\star (c) \bigr)} \begin{cases} > 0 , &
\text{if } x \in \bigl] 0, \gamma^\star (c) \bigr[ , \\
= 0 , & \text{if } x = \gamma^\star (c) , \\ < 0 , & \text{if }
x \in \bigl] \gamma^\star (c),\beta^\star (c) \bigr[ ,
\end{cases} \quad \text{and} \quad
F \bigl( \gamma^\star (c) , \beta^\star (c) \bigr)
 = - c \label{w'-k'}
\een
for all $c \in \mbox{} ]0, c^\star[$.
There exist no other points $0 < \gamma <  \beta <
\infty$ satisfying the system of equations (\ref{BC}).
The functions $\beta^\star$ and $\gamma^\star$
are such that
\begin{gather}
\lim _{c \downarrow 0} \beta^\star (c) =
\lim _{c \downarrow 0} \gamma^\star (c) =
\xlb , \label{beta-gamma-c=0} \\
\lim _{c \uparrow c^\star} \beta^\star (c) = \xub
\quad \text{and} \quad
\lim _{c \uparrow c^\star} \gamma^\star (c)
\begin{cases} > 0 , & \text{if } \lim _{x \downarrow 0}
\frac{R'_\Theta (x)}{\psi' (x)} > 0 \ (\xub = \infty)
, \\ = 0 , & \text{if } \lim _{x \downarrow 0}
\frac{R'_\Theta (x)}{\psi' (x)} = 0 \ (\xub = \infty)
, \\ = 0 , & \text{if } \lim _{x \downarrow 0}
\frac{R'_\Theta (x)}{\psi' (x)} < 0 \ (\xub < \infty) .
\end{cases} \label{lambdac*}
\end{gather}
Furthermore, $c^\star < \infty$ if and only if
\ben
\text{either} \quad {(\rm I)} \ \xub < \infty \quad
\text{or} \quad {(\rm II)} \ \xub = \infty \text{ and }
\lim _{x \uparrow \infty} \frac{\Theta (x)}{r(x)}
> - \infty . \label{c*finite12}
\een
\end{lem}
{\bf Proof.}
%=======
In view of (\ref{z_*}) and (\ref{z^*}) in Lemma~\ref{lem-z*},
we can see that there exists a point $\gamma \in \mbox{}
]0, \beta[$ such that the first equation in (\ref{BC}) holds
true if and only if $\beta \in \mbox{} ]\xlb, \xub[$, in
which case, $\gamma \in \mbox{} ]0, \xlb[$.
In particular, there exists a unique strictly decreasing
function
$\Gamma: \mbox{} ]\xlb, \xub[ \mbox{} \rightarrow
\mbox{} ]0, \xlb[$ such that
\begin{gather}
\frac{R_\Theta' (x)}{\psi' (x)} - \frac{R_\Theta' (\beta)}
{\psi' (\beta)} \begin{cases} > 0 , & \text{if } x \in \bigl]
0, \Gamma (\beta) \bigr[ , \\ = 0 , & \text{if } x =
\Gamma (\beta) , \\ < 0 , & \text{if } x \in \bigl]
\Gamma (\beta), \beta \bigr[ , \end{cases}
\label{Gamma} \\
\biggl( \frac{R_\Theta'}{\psi'} \biggr)' \!\! \bigl( \Gamma
(\beta) \bigr) < 0 , \quad \biggl( \frac{R_\Theta'}{\psi'}
\biggr)' \!\! (\beta) > 0 , \label{Gamma-ineqs} \\
\lim _{\beta \downarrow \xlb} \Gamma (\beta)
= \xlb \quad \text{and} \quad
\lim _{\beta \uparrow \xub} \Gamma (\beta)
\begin{cases} > 0 , & \text{if } \lim _{x \downarrow 0}
\frac{R'_\Theta (x)}{\psi' (x)} > 0 \ (\xub = \infty)
, \\ = 0 , & \text{if } \lim _{x \downarrow 0}
\frac{R'_\Theta (x)}{\psi' (x)} = 0 \ (\xub = \infty)
, \\ = 0 , & \text{if } \lim _{x \downarrow 0}
\frac{R'_\Theta (x)}{\psi' (x)} < 0 \ (\xub < \infty) .
\end{cases} \label{Gamma-lims}
\end{gather}
It follows that the system of equations (\ref{BC}) has
a unique solution $\gamma < \beta$ if and only if
the equation
\ben
F \bigl( \Gamma (\beta) , \beta \bigr) = -c \label{beta*}
\een
has a unique solution $\beta^\star (c) \in \mbox{}
]\xlb, \xub[$.
Using the first expression in (\ref{F}), the identity in
(\ref{Gamma}), the second of the inequalities in
(\ref{Gamma-ineqs}) and the fact that $\psi$ is strictly
increasing, we calculate
\begin{align}
\frac{\di}{\di \beta} F \bigl( \Gamma (\beta) , \beta \bigr)
& = - \biggl( \frac{R_\Theta'}{\psi'} \biggr)' \!\! (\beta) \psi
(\beta) + \biggl( \frac{R_\Theta'}{\psi'} \biggr)' \!\!
\bigl( \Gamma (\beta) \bigr) \psi \bigl( \Gamma (\beta) \bigr)
\Gamma' (\beta) \nonumber \\
& = - \biggl( \frac{R_\Theta'}{\psi'} \biggr)' \!\! (\beta) \psi
(\beta) \Bigl( \psi (\beta) - \psi \bigl( \Gamma (\beta) \bigr)
\Bigr) < 0 . \nonumber
\end{align}
Combining this result with the fact that
\ben
\lim _{\beta \downarrow \xlb} F \bigl( \Gamma (\beta)
, \beta \bigr) = 0 , \label{F@z_*}
\een
which follows from the first limit in
(\ref{Gamma-lims}), we can see that the equation
$F \bigl( \Gamma (\beta) , \beta \bigr) = -c$  has a
unique solution $\beta^\star (c) \in \mbox{} ]\xlb,
\xub[$ if and only if
\ben
c < - \lim_{\beta \uparrow \xub} F \bigl( \Gamma
(\beta) , \beta \bigr) =: c^\star . \label{c*}
\een
We conclude this part of the analysis by noting that
the points $\beta ^\star (c) \in \mbox{} ]\xlb,
\xub[$ and $\gamma ^\star (c) := \Gamma \bigl(
\beta ^\star (c) \bigr) \in \mbox{} ]0, \xlb[$ provide
the unique solution to the system of equations
(\ref{BC}) if $c \in \mbox{} ]0, c^\star[$, while
the system of equations (\ref{BC}) has no solution
such that $0 < \gamma < \beta < \infty$ if $c \geq
c^\star$.
In particular, the inequalities in (\ref{w'-k'}) follow
from the corresponding ones in (\ref{Gamma}).

The fact that the function $\beta \mapsto F \bigl(
\Gamma (\beta) , \beta \bigr)$ is strictly decreasing,
which we have established above, implies that the
function $c \mapsto \beta^\star (c)$ is strictly
increasing because $\beta^\star (c)$ is the unique
solution to equation (\ref{beta*}) for each $c \in
\mbox{} ]0, c^\star[$.
In turn, this result and the fact that $\Gamma$ is strictly
decreasing imply that the function $\gamma ^\star
= \Gamma \circ \beta ^\star$ is strictly decreasing.
The first limit in (\ref{lambdac*}) follows immediately
from (\ref{c*}).
On the other hand, the second limit in (\ref{lambdac*})
follows immediately from the the first limit in
(\ref{lambdac*}) and the second limit in
(\ref{Gamma-lims}).
Furthermore, the identities in (\ref{beta-gamma-c=0})
follow from the first limit in (\ref{Gamma-lims})
and (\ref{F@z_*}).

To establish the equivalence of the inequality
$c^\star < \infty$ with the condition in (\ref{c*finite12}),
we first use the first expression of $F$ in (\ref{F})
and the definition (\ref{c*}) of $c^\star$ to observe that
\be
c^\star = - \lim _{\beta \uparrow \xub} \Bigl( G_\Theta
(\beta) - G_\Theta \bigr( \Gamma (\beta) \bigr) \Bigr) .
\ee
We next use the second limit in (\ref{Gamma-lims})
as well as Lemmas~\ref{lem-G} and~\ref{lem-z*}.
If $\xub < \infty$, then
\be
c^\star = - G_\Theta (\xub) + \lim _{x \downarrow 0}
G_\Theta (x) \stackrel{(\ref{GFlim0})}{=} - G(\xub)
+ \lim _{x \downarrow 0} \frac{\Theta (x)}{r(x)} < \infty ,
\ee
the inequality following because $\Theta / r$ is strictly
increasing in $]0, \xi[$.
If $\xub = \infty$  and $\lim _{x \downarrow 0}
R_\Theta' (x) / \psi' (x) = 0$, then $\lim _{x \downarrow 0}
\Theta (x) / r(x) > - \infty$ thanks to the first implication
in (\ref{z*-impl}).
In this case,
\be
c^\star \stackrel{(\ref{GFlimoo})}{=} - \lim _{x \uparrow \infty}
\frac{\Theta (x)}{r(x)} + \lim _{x \downarrow 0}
\frac{\Theta (x)}{r(x)} \begin{cases} < \infty , & \text{if }
\lim_{x \uparrow \infty} \frac{\Theta (x)}{r(x)} > - \infty ,
\\ = \infty , & \text{if } \lim_{x \uparrow \infty}
\frac{\Theta (x)}{r(x)} = - \infty , \end{cases}
\ee
where we have also used the assumption that
$\Theta / r$ is strictly decreasing in $]\xi, \infty[$.
Finally, if $\xub = \infty$ and $\lim
_{x \downarrow 0} \frac{R_\Theta' (x)}{\psi' (x)} > 0$,
then $\lim _{x \uparrow \infty} \Gamma (x)> 0$ (see
(\ref{Gamma-lims})),
\be
c^\star \stackrel{(\ref{GFlimoo})}{=} - \lim _{x \uparrow \infty}
\frac{\Theta (x)}{r(x)} + \lim _{x \uparrow \infty}G \bigr(
\Gamma (x) \bigr) \begin{cases} < \infty , & \text{if }
\lim_{x \uparrow \infty} \frac{\Theta (x)}{r(x)} > - \infty ,
\\ = \infty , & \text{if } \lim_{x \uparrow \infty}
\frac{\Theta (x)}{r(x)} = - \infty , \end{cases}
\ee
and the proof is complete.
\mbox{} \hfill $\Box$
\bigskip

%==============================================
In light of (\ref{OS2}), (\ref{OS3}) and the previous
lemma, we now establish the following result, which
provides the solution to the HJB equation (\ref{HJB})
identifying with the control problem's value function
when a $\beta$-$\gamma$ strategy with $\gamma <
\beta$ in $]0, \infty[$ is indeed optimal.

\begin{lem} \label{lem-existGamma}
Consider the stochastic control problem formulated in
Section~\ref{sec:pr-form} and suppose that the
point $\xlb$ introduced in Lemma~\ref{lem-z*}.(i) is
finite.
Also, fix any $c \in \mbox{} ]0,c^\star[$, where
$c^\star>0$ is as in Lemma~\ref{lem:FB}.
The function $w$ defined by
\ben
w(x) = \begin{cases} R_h (x) + \Bigl( K_\infty -
\frac{R_\Theta' (\beta^\star)}{\psi' (\beta^\star)}
\Bigr) \psi(x) , & \text{for } x \in \mbox{} ]0, \beta^\star[
, \\ w(\gamma^\star) + \int _{\gamma^\star}^x k(s)
\, \di s - c , & \text{for } x \in [\beta^\star, \infty[ ,
\end{cases} \label{w-sol}
\een
where we write $\gamma^\star$ and $\beta^\star$
in place of the points $\gamma^\star (c)$ and
$\beta^\star (c)$ given by Lemma~\ref{lem:FB},
is $C^1$ in $]0, \infty[$ and $C^2$ in $]0, \infty[
\mbox{} \setminus \{ \beta^\star \}$.
Furthermore, this function is a solution to the HJB
equation (\ref{HJB}) that is bounded from below.
\end{lem}
{\bf Proof.}
%=======
The boundedness from below of $w$ follows immediately
from Assumption~\ref{A3}, the conditions in (i) and (ii) of Assumption~\ref{A5}
and Lemma~\ref{lem:R_F}.(i).

By construction, we will establish all of the lemma's
other claims if we prove that
\begin{align}
- c + \int _u^x \bigl( k(s) - w'(s) \bigr) \, \di s \leq 0 & \quad
\text{for all } 0 < u < x < \beta^\star \label{HJB-ineq1} \\
\text{and} \quad \lscr w(x) + h(x) \leq 0 & \quad
\text{for all } x > \beta^\star . \label{HJB-ineq2}
\end{align}
To this end, we use the first expression of $w$
in (\ref{OS3}) and (\ref{w'-k'}) to note that
\be
k(s) - w'(s) = \psi'(s) \biggl(
\frac{R_\Theta' (\beta^\star)}{\psi' (\beta^\star)}
- \frac{R_\Theta' (s)}{\psi' (s)} \biggr)
\begin{cases} < 0 & \text{for all } s \in \mbox{} ]0,
\gamma^\star[ , \\ > 0 & \text{for all } s \in \mbox{}
]\gamma^\star, \beta^\star[ .  \end{cases}
\ee
The inequality (\ref{HJB-ineq1}) follows from this
observation and the fact that
\be
- c + \int _{\gamma^\star}^{\beta^\star} \bigl( k(s)
- w'(s) \bigr) \, \di s = 0 .
\ee

To show (\ref{HJB-ineq2}), we first use the expression
\be
w(x) = w(\beta^\star) + \int _{\beta^\star}^x k(s) \, \di s ,
\quad \text{for } x > \beta^\star ,
\ee
the definition (\ref{Theta}) of $\Theta$ in
Assumption~\ref{A5} and the first expression in
(\ref{OS3}) to calculate
\begin{align}
\lscr w(x) + h(x) & = - r(x) w(\beta^\star) + \lscr
\biggl( \int_0^{\boldsymbol{\cdot}}
k(s) \, \di s \biggr) (x) + r(x) \int _0^{\beta^\star} k(s)
\, \di s + h(x) \nonumber \\
& = \Theta (x) - r( x)\biggl( R_\Theta (\beta^\star)
- \frac{R_\Theta' (\beta^\star)}{\psi' (\beta^\star)}
\psi (\beta^\star) \biggr) \nonumber \\
& = \Theta (x) -r(x) G_\Theta (\beta^\star) ,
\label{Theta-rG}
\end{align}
where $G_\Theta$ is given by (\ref{G}) in
Lemma~\ref{lem-G} for $F=\Theta$.
In view of the calculations
\be
G_\Theta' (x) = - \psi (x) \frac{\di}{\di x}
\frac{R_\Theta' (x)}{\psi' (x)} = - \frac{2C r(x) p'(x) \psi (x)}
{( \sigma (x) \psi' (x) )^2} \biggl( \int _0^x \Theta
(s) \Psi (s) \, \di s - \frac{\Theta (x)}{r(x)}
\frac{\psi' (x)}{C p'(x)} \biggr) ,
\ee
the inequalities (\ref{z_*}) in Lemma~\ref{lem-z*} and
the fact that $\beta^\star > \xlb$, we can see that
\ben
G_\Theta (x) < G_\Theta (\beta^\star) \quad
\text{and} \quad \int _0^x \Theta (s) \Psi (s) \, \di s >
\frac{\Theta (x)}{r(x)} \frac{\psi' (x)}{C p'(x)} \quad
\text{for all } x > \beta^\star .  \label{HJB-ineq22}
\een
The second of these inequalities and the second expression
of $G_\Theta$ in (\ref{G}) imply that
\be
G_\Theta (x) = \frac{C p'(x)}{\psi' (x)} \int _0^x \Theta (s)
\Psi (s) \, \di s > \frac{\Theta (x)}{r(x)} \quad \text{for all }
x > \beta^\star .
\ee
However, this result, (\ref{Theta-rG}) and the first inequality
in (\ref{HJB-ineq22}) yield
\be
\lscr w(x) + h(x) < r(x) \biggl( \frac{\Theta (x)}{r(x)} -
G_\Theta (x) \biggr) < 0 \quad \text{for all } x > \beta^\star ,
\ee
and (\ref{HJB-ineq2}) follows.
\mbox{} \hfill $\Box$
\bigskip

%==============================================
To proceed further, we assume that the problem
data is such that $c^\star < \infty$, which is the case if and
only if one of the two conditions of (\ref{c*finite12}) in
Lemma~\ref{lem:FB} holds true.
In the first case, when $\xub < \infty$, the limits in
(\ref{lambdac*}) suggest the possibility for the function
$w$ defined by (\ref{OS2}) and (\ref{OS3}) for $\gamma
= 0$ and some $\beta > \xub$ to provide a solution
to the HJB equation (\ref{HJB}) that identifies with the
control problem's value function.
In this case, a free-boundary condition such as
(\ref{condg}) is not relevant anymore and we are faced
with only the free-boundary condition (\ref{wgabe}) with
$\gamma = 0$, which is equivalent to the equation
$F(0, \beta) = -c$, where $F$ is defined by (\ref{F}).

\begin{lem} \label{lem:beta-0}
Consider the stochastic control problem formulated in
Section~\ref{sec:pr-form} and suppose that the
point $\xlb$ introduced in Lemma~\ref{lem-z*}.(i) is
finite.
Also, suppose that the problem
data is such that $\xub < \infty$, where $\xub$ is
defined by (\ref{z^*}) in Lemma~\ref{lem-z*}.
The following statements hold true:
\smallskip

\noindent {\rm (I)}
There exists $c^\circ \in \mbox{} ]c^\star, \infty]$ and
a strictly increasing function $\beta^\circ:  [c^\star, c^\circ[
\mbox{} \rightarrow [\xub, \infty[$ such that  
\ben
F \bigl( 0, \beta^\circ (c) \bigr) = - c \text{ for all } c \in
[c^\star, c^\circ[ \quad \text{and} \quad
\lim _{c \uparrow c^\circ} \beta^\circ (c) = \infty ,
\label{beta^o}
\een
where $c^\star \in \mbox{} ]0, \infty[$ is as in 
Lemma~\ref{lem:FB}.
\smallskip

\noindent {\rm (II)}
$c^\circ = \infty$ if and only if
$\lim _{x \uparrow \infty} \Theta(x) / r(x) = -\infty$.
\smallskip

\noindent {\rm (III)}
Given any $c \in [c^\star, c^\circ[$, the function $w$
defined by
\ben
w(x) = \begin{cases} R_h (x) + \Bigl( K_\infty -
\frac{R_\Theta' (\beta^\circ)}{\psi' (\beta^\circ)}
\Bigr) \psi(x) , & \text{for } x \in \mbox{} ]0, \beta^\circ[
, \\ R_h (0) + \Bigl( K_\infty -
\frac{R_\Theta' (\beta^\circ)}{\psi' (\beta^\circ)} \Bigr)
\psi (0) + \int _0^x k(s) \, \di s - c , & \text{for } x \in
[\beta^\circ, \infty[ , \end{cases} \label{w-sol-0}
\een
where we write $\beta^\circ$ in place of $\beta^\circ (c)$,
is $C^1$ in $]0, \infty[$ and $C^2$ in $]0, \infty[ \mbox{}
\setminus \{ \beta^\circ \}$.
Furthermore, this function is a solution to the HJB
equation (\ref{HJB}) that is bounded from below.
\end{lem}
{\bf Proof.} 
%=======
The definition of $G_\Theta$ as in (\ref{G}), the limits
(\ref{GFlim0}) in Lemma~\ref{lem-G} and the implications
(\ref{z*-impl}) in Lemma~\ref{lem-z*} imply that the
limit $\lim _{x \downarrow 0} G_\Theta (x)$ exists in
$\bbr$ thanks to Assumption~\ref{A5}.(iii).
On the other hand, (\ref{z_*}) and (\ref{z*-equiv}) in
Lemma~\ref{lem-z*} imply that the limit $\lim
_{x \downarrow 0} R_\Theta' (x) / \psi' (x)$ exists
in $]{-\infty}, 0[$.
In view of these observations and the definition of
$G_\Theta$ as in (\ref{G}), we can see that the limit
$\lim _{x \downarrow 0} R_\Theta (x)$ exists in
$\bbr$.
Therefore, the limit $R_h (0) := \lim _{x \downarrow 0}
R_h (x)$ exists in $\bbr$ thanks to (\ref{Rh-RTheta}).
It follows that the function $w$ is well-defined.

The second expression in (\ref{w'-k'}) and the limits
in (\ref{lambdac*}) imply that
\be
F(0, \xub) \equiv G_\Theta (\xub) - \lim _{x \downarrow 0}
G_\Theta (x) = - c^\star \in \mbox{} ]{-\infty}, 0[ .
\ee
Part~(I) of the lemma follows from this observation and
the calculation
\be
\frac{\di}{\di \beta} F (0, \beta) = - \biggl(
\frac{R_\Theta'}{\psi'} \biggr)' \!\! (\beta) \bigl( \psi (\beta)
- \psi (0) \bigr) < 0 \quad \text{for all } \beta \geq \xub ,
\ee
where the inequality follows from (\ref{z_*}) in
Lemma~\ref{lem-z*} and the fact that the strictly positive
function $\psi$ is strictly increasing, for $c^\circ = - \lim
_{\beta \uparrow \infty} F(0, \beta)$.
Furthermore, this definition of $c^\circ$,
Assumption~\ref{A5}.(iii) and the limits (\ref{GFlimoo}) in
Lemma~\ref{lem-G} imply immediately part~(II) of the
lemma.

Finally, we can show the rest of the claims on $w$ by
using exactly the same  arguments as in the proof
of Lemma~\ref{lem-existGamma} (see (\ref{HJB-ineq1})
and (\ref{HJB-ineq2}) in particular).
\mbox{}\hfill $\Box$
\bigskip

%==============================================
To close the ``gap'' in the parameter space, we still
need to derive a solution to the HJB equation (\ref{HJB})
if
\be
\xlb < \xub = \infty, \ c^\star < \infty \text{ and }
c \geq c^\star , \quad \text{or} \quad
\xub < \infty, \ c^\circ < \infty \text{ and }
c \geq c^\circ , \quad \text{or} \quad
\xlb = \infty \text{ and } c > 0 .
\ee
In the first case, the first limit in (\ref{lambdac*}) implies
that $\lim _{c \uparrow c^\star} \beta^\star (c) = \infty$.
In the second case, the limit in (\ref{beta^o}) implies
that $\lim _{c \uparrow c^\circ} \beta^\circ (c) = \infty$.
In all cases, we are faced with the possibility for the
problem's value function to identify with a solution
to the ODE $\lscr w(x) + h(x)= 0$ for all $x > 0$.

\begin{lem} \label{lem:oo-0}
Consider the stochastic control problem formulated in
Section~\ref{sec:pr-form} and suppose that the
problem data is such that one of the following cases
holds true:
\smallskip

\noindent {\rm (a)}
The point $\xlb$ introduced in Lemma~\ref{lem-z*}.(i) is
finite,
\ben
\lim _{x \uparrow \infty} \frac{\Theta (x)}{r(x)} > - \infty
\quad \Leftrightarrow \quad
\text{either } (\xub = \infty \text{ and } c^\star < \infty)
\text{ or } (\xub < \infty \text{ and } c^\circ < \infty)
\label{final-equiv}
\een
and $c \geq c^\star$ or $c \geq c^\circ$, depending
on the case in (\ref{final-equiv}).
\smallskip

\noindent {\rm (b)}
The point $\xlb$ introduced in Lemma~\ref{lem-z*}.(i) is
equal to infinity.
\smallskip

\noindent
In either of these two cases, the function $w$ defined
by
\ben
w(x) = R_h (x) + K_\infty \psi(x) , \quad \text{for } x>0 ,
\label{w-final}
\een
is a $C^2$ solution to the HJB equation (\ref{HJB}) that
is bounded from below.
\end{lem}
{\bf Proof.} 
%=======
The equivalence (\ref{final-equiv}) follows immediately
from the statement related to (\ref{c*finite12}) in
Lemma~\ref{lem:FB} and part~(II) of
Lemma~\ref{lem:beta-0}.
On the other hand, the boundedness from below of $w$
follows immediately from Assumption~\ref{A3}, the
conditions in (i) and (ii) of Assumption~\ref{A5} and
Lemma~\ref{lem:R_F}.(i).

To establish the fact that $w$ satisfies the HJB equation
(\ref{HJB}), we have to show that
\ben
\int _u^x \bigr( k(s) - w'(s) \bigr) \, \di s \leq c \quad
\Leftrightarrow \quad R_\Theta (u) - R_\Theta (x)
\leq c \quad\text{for all } 0 < u < x < \infty ,
\label{HJB-ineq-f}
\een
where the equivalence follows from the identity
(\ref{Rh-RTheta}) and the definition (\ref{w-final})
of $w$.
To this end, fix any $u<x$ in $]0, \infty[$.
First, suppose that $\xub = \infty$ and $c^\star < \infty$.
In this case, the limits in (\ref{lambdac*}) imply that
$x < \beta^\star (c)$ for all $c<c^\star$ sufficiently
close to $c^\star$.
For such a $c$, the identity (\ref{Rh-RTheta}) and
the fact that the function $w$ defined by (\ref{w-sol})
in Lemma~\ref{lem-existGamma} satisfies the HJB
equation (\ref{HJB}) imply that
\be
R_\Theta (u) - R_\Theta (x) \leq c +
\frac{R_\Theta' \bigl( \beta^\star (c) \bigr)}
{\psi' \bigl( \beta^\star (c) \bigr)} \bigl( \psi (u)
- \psi (x) \bigr) .
\ee
Passing to the limit as $c \uparrow c^\star$ and using
the fact that $\lim _{c \uparrow c^\star} \beta^\star (c)
= \infty$ together with the limit in (\ref{z_*}), we can see
that $R_\Theta (u) - R_\Theta (x) \leq c^\star$.
It follows that (\ref{HJB-ineq-f}) holds true for all
$c \geq c^\star$.

If $\xub < \infty$, $c^\circ < \infty$ and $c \geq c^\circ$,
then we can show that the function $w$ given by
(\ref{w-final}) satisfies the HJB equation (\ref{HJB})
in exactly the same way using the results of
Lemma~\ref{lem:beta-0}.

Finally, suppose that the point $\xlb$ introduced in
Lemma~\ref{lem-z*}.(i) is equal to infinity and
consider any points $u < x < \beta$ in $]0, \infty[$.
In this case, the inequalities in (\ref{z_*}) imply that
\be
R_\Theta' (s) - \frac{R_\Theta' (\beta)}{\psi' (\beta)}
\psi' (s) = \psi' (s) \biggl( \frac{R_\Theta' (s)}{\psi' (s)}
- \frac{R_\Theta' (\beta)}{\psi' (\beta)} \biggr) > 0
\quad \text{for all } s < \beta .
\ee
In view of this observation, we can see that
\be
R_\Theta (u) - R_\Theta (x) \leq
- \frac{R_\Theta' (\beta)}{\psi' (\beta)} \bigl( \psi (x)
- \psi (u) \bigr) .
\ee
Passing to the limit as $\beta \rightarrow \infty$, we
can see that $R_\Theta (u) - R_\Theta (x) \leq 0$,
thanks to the limit in (\ref{z_*}).
It follows that (\ref{HJB-ineq-f}) holds true for all
$c>0$.
\mbox{}\hfill $\Box$
\bigskip

%==============================================
We conclude the section with the main result of the paper. 

\begin{thm} \label{thm:SolPb-lem}
Consider the stochastic control problem formulated in
Section~\ref{sec:pr-form}.
Depending on the problem data, the function $w$ defined
by (\ref{w-sol}), (\ref{w-sol-0}) or (\ref{w-final}) in
Lemmas~\ref{lem-existGamma}, \ref{lem:beta-0}
or~\ref{lem:oo-0}, respectively, identifies with the control
problem's value function, namely,
\ben
w(x) = \sup_{\zeta \in \acal} J_x(\zeta) . \label{supJ=w}
\een
Furthermore, the following cases hold true:
\smallskip

\noindent {\rm (I)}
If the problem data is as in Lemma~\ref{lem-existGamma},
then the $\beta$-$\gamma$ strategy characterised
by the points $\beta^\star$ and $\gamma^\star$ in
Lemma~\ref{lem-existGamma} is optimal.
\smallskip

\noindent {\rm (II)}
If the problem data is as in Lemma~\ref{lem:beta-0},
then there exists no optimal strategy.
In this case, if $(\varepsilon_n)$ is any sequence
such that $\varepsilon_1 < \beta^\circ$ and $\lim
_{n \uparrow \infty} \varepsilon_n = 0$, then
the $\beta$-$\gamma$ strategies characterised
by the points $\beta = \beta^\circ$ and $\gamma
= \varepsilon_n$, where $\beta^\circ$ is as in
Lemma~\ref{lem:beta-0}, provide a sequence of
$\varepsilon$-optimal strategies.
\smallskip

\noindent {\rm (III)}
If the problem data is as in Lemma~\ref{lem:oo-0}
and $K_\infty = 0$, then $\zeta^\star = 0$ is an
optimal strategy.
\smallskip

\noindent {\rm (IV)}
If the problem data is as in Lemma~\ref{lem:oo-0}
and $K_\infty > 0$, then there exists no optimal strategy.
In this case, if $\gamma$ is an arbitrary point in $]0, \infty[$
and $(\varepsilon_n)$ is any sequence such that
$\varepsilon _1^{-1} > \gamma$ and $\lim
_{n \uparrow \infty} \varepsilon _n^{-1} = \infty$, then
the $\beta$-$\gamma$ strategies characterised by the
points $\beta = \varepsilon_n^{-1}$ and $\gamma$
provide a sequence of $\varepsilon$-optimal strategies.
\end{thm}
{\bf Proof.}
%=======
Fix any initial value $x>0$, consider any admissible controlled
process $\zeta \in \acal$ and denote by $X^\zeta$ the
associated solution to the SDE (\ref{SDE}).
Using It\^{o}'s formula,  we obtain
\begin{align}
e^{-\Lambda _T^\zeta} w(X _T^\zeta) = w(x) + \int _0^T
e^{-\Lambda _t^\zeta} \lscr w(X_t^\zeta) \, \di t + \sum
_{0 \leq t \leq T} e^{-\Lambda _t^\zeta} \bigl(
w(X_t^\zeta) - w(X_{t-}^\zeta) \bigr) {\bf 1}
_{\{ \Delta \zeta_t > 0 \}} + M_T^\zeta , \nonumber
\end{align}
where
\be
M_T^\zeta = \int _0^T e^{-\Lambda _t^\zeta} \sigma
(X_t^\zeta) w' (X_t^\zeta) \, \di W_t .
\ee
Since $\Delta X_t^\zeta \equiv X_t^\zeta - X^\zeta_{t-}
= - \Delta \zeta_t \leq 0$, we can see that
\begin{align}
w (X_t^\zeta) - w (X_{t-}^\zeta) + \int _0^{\Delta \zeta_t}
& k(X_{t-}^\zeta - u) \, \di u - c {\bf 1} _{\{ \Delta \zeta_t > 0 \}}
\nonumber \\
& = \biggl( \int _{X_{t-}^\zeta - \Delta \zeta_t}^{X_{t-}^\zeta}
\bigl( k(u) - w' (u) \bigr) \, \di u - c \biggr) {\bf 1}
_{\{ \Delta \zeta_t > 0 \}} . \nonumber
\end{align}
In view of these observations and the fact that $w$
satisfies the HJB equation (\ref{HJB}), we derive
\begin{align}
\int _0^T e^{-\Lambda _t^\zeta} h(X_t^\zeta) & \, \di t
+ \sum _{t \in [0,T]} e^{-\Lambda _t^\zeta} \biggl( \int
_0^{\Delta \zeta_t} k(X_{t-}^\zeta -u) \, \di u - c {\bf 1}
_{\{ \Delta \zeta_t > 0 \}} \biggr) \nonumber \\
= \mbox{} & w(x) - e^{-\Lambda _T^\zeta}
w(X_T^\zeta) + \int _0^T e^{-\Lambda _t^\zeta} \Bigl(
\lscr w(X_t^\zeta) + h(X_t^\zeta) \Bigr) \, \di t
\nonumber \\
& + \sum _{0 \leq t \leq T} \biggl( e^{-\Lambda _t^\zeta}
\int _{X_{t-}^\zeta - \Delta \zeta_t}^{X_{t-}^\zeta} \bigl( k(u) -
w' (u) \bigr) \, \di u -c \biggr) {\bf 1} _{\{ \Delta \zeta_t > 0 \}}
+ M_T^\zeta \nonumber \\
\leq \mbox{} & w(x) - e^{-\Lambda _T^\zeta}
w(X_T^\zeta) + M_T ^\zeta. \label{VTineq1}
\end{align}

We next consider any sequence $(\tau_n)$ of bounded
localising times for the local martingale $M^\zeta$.
Recalling Assumption~\ref{A5}.(ii) as well as the fact
that $h$ and $w$ are both bounded from below,
we use Fatou's lemma, the monotone convergence theorem
and the admissibility condition (\ref{AC}) to observe that
(\ref{VTineq1}) implies that
\begin{align}
J_x (\zeta) & \leq \liminf _{n \uparrow \infty} \bbe_x \Biggl[
\int _0^{\tau_n} e^{-\Lambda _t^\zeta} h (X_t^\zeta) \,
\di t + \sum _{t \in [0, \tau_n]} e^{-\Lambda _t^\zeta}
\biggl( \int _0^{\Delta \zeta_t} k(X_{t-}^\zeta -u) \, \di u
- c {\bf 1} _{\{ \Delta \zeta_t > 0 \}} \biggr) \Biggr]
\nonumber \\
& \leq \lim _{n \uparrow \infty} \bbe_x \Bigl[ w(x) +
e^{-\Lambda _{\tau_n}^\zeta} w^- (X_{\tau_n}^\zeta)
\Bigr] = w(x) , \label{J<w}
\end{align}
where $w^- (x) = - \min \bigl\{ 0, w(x) \bigr\}$.
\smallskip

{\em Proof of (I)\/.}
First, consider any $x \in \mbox{} ]0, \beta[$.
In view of the results in Theorem~\ref{lem:beta-gamma},
the $\beta$-$\gamma$ strategy $\zeta^\star$ characterised
by the points $\beta^\star$ and $\gamma^\star$ is such
that
\ben
J_x (\zeta^\star) = R_h (x) +
\frac{\psi (x)}{\psi (\beta^\star) - \psi (\gamma^\star)}
\biggl( R_h (\gamma^\star) - R_h (\beta^\star) + \int
_{\gamma^\star}^{\beta^\star} k(s) \, \di s - c \biggr) .
\label{J*}
\een
On the other hand, the identity $F(\gamma^\star, \beta^\star)
= - c$ and the definition (\ref{F}) of $F$ imply that
\be
\frac{R_\Theta' (\beta^\star)}{\psi' (\beta^\star)} =
\frac{R_\Theta (\beta^\star) - R_\Theta (\gamma^\star) + c}
{\psi (\beta^\star) - \psi (\gamma^\star)} .
\ee
In view of the identity (\ref{Rh-RTheta}), this expression
is equivalent to
\be
\frac{R_\Theta' (\beta^\star)}{\psi' (\beta^\star)} =
\frac{1}{\psi (\beta^\star) - \psi (\gamma^\star)}
\biggl( R_h (\beta^\star) - R_h (\gamma^\star)
- \int _{\gamma^\star}^{\beta^\star} k(s) \, \di s + c
+ K_\infty \bigl( \psi (\beta^\star) - \psi (\gamma^\star)
\bigr) \biggr) .
\ee
However, this result, the definition (\ref{w-sol}) of $w$
and (\ref{J*}) imply that $J_x (\zeta^\star) = w(x)$,
which, combined with (\ref{J<w}), establishes
(\ref{supJ=w}) as well as the optimality of $\zeta^\star$.
The corresponding claims for $x \geq \beta$ are
immediate.
\smallskip

{\em Proof of (II)\/.}
In this case, the identity $F(0, \beta^\circ) = -c$
implies that the sequence $(c_n)$ defined by $c_n
= -F(\varepsilon_n, \beta^\circ)$ is such that
$\lim _{n \uparrow \infty} c_n = c$.
By following reasoning similar to the one in the
previous part of the proof, we can see that, given any
$x \in \mbox{} ]0, \beta[$, the $\beta$-$\gamma$ strategy
$\zeta^{\varepsilon_n}$ characterised by the points
$\beta = \beta^\circ$ and $\gamma = \varepsilon_n$
is such
that
\be
J_x (\zeta^{\varepsilon_n}) = w(x) -
\frac{(c - c_n) \psi (x)}{\psi (\beta^\circ) - \psi (\varepsilon_n)} ,
\ee
and the required results follow.
\smallskip

{\em Proof of (III)\/.}
This case follows immediately from (\ref{J<w}) and
the probabilistic expression of $R_h$ as in (\ref{RF}).
\smallskip

{\em Proof of (IV)\/.}
In view of the results in Theorem~\ref{lem:beta-gamma},
the $\beta$-$\gamma$ strategy $\zeta^{\varepsilon_n}$
characterised by the points $\beta = \varepsilon_n^{-1}$
and $\gamma$ is such that
\be
J_x (\zeta^{\varepsilon_n}) = R_h (x) +
\frac{\psi (x)}{\psi (\varepsilon_n^{-1}) - \psi (\gamma)}
\biggl( R_h (\gamma) - R_h (\varepsilon_n^{-1})
+ \int _\gamma^{\varepsilon_n^{-1}} k(s)
\, \di s - c \biggr) .
\ee
Combining this observation with the second limit in
(\ref{limitRF/phiRF/psi}) and the definition (\ref{Koo})
of $K_\infty$, we can see that $\lim _{n \uparrow \infty}
J_x (\zeta^{\varepsilon_n}) = R_h (x) + K_\infty \psi (x)$.
However, this limit and (\ref{J<w}) imply the required
results.
\mbox{}\hfill $\Box$

\begin{rem} \label{rem:nat/entr} {\rm
Suppose that we enlarged the family of admissible strategies
to allow for switching the system off.
In particular, suppose that we allowed for the controlled
process $X^\zeta$ to hit 0 at some time and be absorbed
by 0 after that time.
In this context, we would face the HJB equation
\begin{align}
\max \Biggl\{ \lscr w(x) + h(x) , \ -c + \sup _{z \in [0,x[}
\int _{x-z}^x \bigl( k(s) - w'(s) \bigr) \, \di s , & \nonumber \\
-w(0) - c + \frac{h(0)}{r(0)} + \int _0^x \bigl( k(s) - w'(s) \bigr)
\, \di s & \Biggr\} = 0 , \label{HJBs}
\end{align}
where we assume that both of the limits $h(0) :=
\lim _{x \downarrow 0} h(x)$ and $r(0) := \lim
_{x \downarrow 0} r(x)$ exist in $\bbr$, instead of
just the limit $\lim _{x \downarrow 0} h(x) / r(x)$.
The third term of this HJB equation incorporates the
inequality
\be
w(x) \geq -c + \int _0^x k(s) \, \di s + \int _0^\infty
e^{-r(0) s} h(0) \, \di s
\ee
that should hold with equality for those values $x$
of the state space at which it is optimal to switch off the
system.

In view of the second limit in (\ref{0-natural}) and
Lemma~\ref{lem:R_F}.(ii), if 0 is a natural boundary
point, then, in all of the cases appearing in
Lemmas~\ref{lem-existGamma}-\ref{lem:oo-0},
\be
w(0) = R_h (0) = \frac{h(0)}{r(0)}
\ee
and the inequality associated with the third term of
(\ref{HJBs}) follows from the one associated with the
second term of (\ref{HJBs}).
In view of this observation, we can see that the
results of Theorem~\ref{thm:SolPb-lem} hold true
with the following modification:
in Case~II, the $\beta$-0 strategy that switches off
the system as soon as the uncontrolled process $X$
takes any value greater than or equal to $\beta =
\beta^\circ$ is optimal.
In Case~IV of the theorem, an optimal strategy still
does not exist.

The situation is entirely different if 0 is an entrance
boundary point.
In this case, Theorem~\ref{thm:SolPb-lem} with a
modification such as the one in the previous paragraph
still provides the solution to the control problem if the
problem data is such that the solution $w$ to the HJB
equation (\ref{HJB}) satisfies the inequality
$w(0) \geq h(0) / r(0)$.
In Example~\ref{ex8} in the next section, we can see
that this inequality may or may not be true.
In particular, a $\beta$-0 strategy that may switch the
system off can indeed be optimal and be associated
with a payoff that is strictly greater than the value
function derived in Theorem~\ref{thm:SolPb-lem}.
Investigating the solution to the control problem if
we allowed for the system to be switched off would
require substantial extra analysis that goes beyond
the scope of the present article.
\mbox{} \hfill $\Box$ }
\end{rem}

%%%%%%%%%%%%%%%%%%%%%%%%%%%%%%%%%%
\section{Examples}
\label{sec:ex}

The first four examples that we consider in this section
present choices for the problem data that satisfy our
assumptions.
In these examples, the functions $r$ and $k$ are strictly
positive constants, so the function $\Theta$ introduced
in Assumption~\ref{A5} takes the form
\be
\Theta (x) = h(x) + k b(x) - rk x .
\ee
Furthermore, $\lim _{x \uparrow \infty} \Theta (x) = -\infty$
in each of the Examples~\ref{ex1}-\ref{ex4}, which implies
that $\xlb < \infty$ thanks to Lemma~\ref{lem-z*}.(ii),
where $\xlb \in \mbox{} ]\xi, \infty]$ is as in
Lemma~\ref{lem-z*}.(i).

\begin{ex} \label{ex1} {\rm
Suppose that the uncontrolled dynamics of the state process
are modelled by the SDE
\ben
\di X_t = bX_t \, \di t + \sigma X_t \, \di W_t , \quad X_0 =
x>0 , \label{gBm}
\een
for some constants $b$ and $\sigma > 0$.
Furthermore, if $r>b$ and $h$ is any bounded from below
strictly concave function such that
\be
\lim _{x \downarrow 0} h'(x) > k (r-b)
\quad \text{and} \quad
\lim _{x \uparrow \infty} h'(x) = 0 ,
\ee
then $\Theta$ is strictly concave and satisfies the
requirements of Assumption~\ref{A5}.
} \end{ex}

\begin{ex} \label{ex2} {\rm
Suppose that the uncontrolled dynamics of the state process
are modelled by the SDE
\be
\di X_t = \upkappa (\upgamma - X_t) X_t \, \di t
+ \upsigma X_t^\ell \, \di W_t , \quad X_0 = x  > 0 ,
\ee
for some strictly positive constants $\upkappa$, $\upgamma$,
$\upsigma$ and $\ell \in [1, \frac{3}{2}]$.
Note that the celebrated stochastic Verhulst-Pearl logistic
model of population growth arises in the special case
$\ell = 1$.
Assumptions~\ref{A1}-\ref{A3} hold true
if $\ell \in \mbox{} ]1, \frac{3}{2}]$ or if $\ell = 1$
and $k \upgamma - \frac{1}{2} \sigma^2 > 0$.
Furthermore, if $h$ is any bounded from below
concave function such that
\be
\lim _{x \downarrow 0} h'(x) > k (r - \kappa \gamma) ,
\ee
then $\Theta$ is strictly concave and satisfies the
requirements of Assumption~\ref{A5}.
\mbox{} \hfill $\Box$ }
\end{ex}

\begin{ex} \label{ex3} {\rm
Suppose that the uncontrolled dynamics of the state process
are modelled by the SDE
\be
\di X_t = \biggl( \upkappa \upgamma + \frac{1}{2}
\sigma ^2 - \upkappa \ln (X_t) \biggr) X_t \, \di t +
\upsigma X_t \, \di W_t , \quad X_0 = x  > 0 ,
\ee
for some constants $\upkappa, \upgamma, \upsigma
> 0$, namely, the logarithm of the uncontrolled state
process is the Ornstein-Uhlenbeck process given by
\be
\di \ln (X_t) = \upkappa \bigl( \upgamma - \ln (X_t) \bigr)
\, \di t + \upsigma \, \di W_t , \quad \ln (X_0) = \ln (x)
\in \bbr .
\ee
Furthermore, if $h$ is any bounded from below
concave function, then $\Theta$ is strictly concave
and satisfies the requirements of Assumption~\ref{A5}.
\mbox{} \hfill $\Box$ }
\end{ex}

\begin{ex} \label{ex4} {\rm
Suppose that the uncontrolled dynamics of the state process
are modelled by the SDE
\be
\di X_t = \upkappa (\upgamma - X_t) \, \di t + \upsigma
X_t^\ell \, \di W_t , \quad X_0 = x  > 0 ,
\ee
for some strictly positive constants $\upkappa$, $\upgamma$,
$\upsigma$ and $\ell \in [ \frac{1}{2} , 1]$.
Note that, in the special case that arises for $\ell = \frac{1}{2}$
and $\upkappa \upgamma - \frac{1}{2} \upsigma^2 > 0$,
the process $X$ identifies with the short rate
process in the Cox-Ingersoll-Ross interest rate model.
Assumptions~\ref{A1}-\ref{A3} hold true
if $\ell \in \mbox{} ]\frac{1}{2} , 1]$ or if $\ell = \frac{1}{2}$
and $k\upgamma - \frac{1}{2} \sigma^2 > 0$.
Furthermore, if $h$ is any bounded from below
strictly concave function such that
\be
\lim _{x \downarrow 0} h'(x) > k (r + \kappa)
\quad \text{and} \quad
\lim _{x \uparrow \infty} h'(x) = 0 ,
\ee
then $\Theta$ is strictly concave and
satisfies the requirements of Assumption~\ref{A5}.
\mbox{} \hfill $\Box$ }
\end{ex}

The next three examples illustrate the four different cases
that appear in Theorem~\ref{thm:SolPb-lem}, our main result.
In the next three ones, $X$ is the geometric Brownian motion
that is given by (\ref{gBm}).
In this context, it is well-known that
\be
\varphi (x) = x^m, \quad \psi (x) = x^n \quad \text{and}
\quad p'(x) = x^{m+n-1} ,
\ee
where the constants $m<0<n$ are given by
\be
m, n = \frac{1}{2} - \frac{b}{\sigma^2} \mp
\sqrt{\left(\frac{1}{2}-  \frac{b}{\sigma^2}\right)^2
+ \frac{2r}{\sigma^2} } ,
\ee
while the constant $C$ defined by (\ref{C}) is equal to
$n-m$.
Furthermore, the identities
\be
mn = -\frac{2r}{\sigma^2} \quad \text{and} \quad m+n =
1 - \frac{2b}{\sigma^2}
\ee
hold true, while
\be
r<b \ \Leftrightarrow \ 0<n<1 \quad \text{and} \quad
b=r \ \Leftrightarrow \ 1=n .
\ee

%----------------------------------------------------------------------
\begin{ex} \label{ex5} {\rm
Suppose that $r > b$ and consider the functions
\be
h(x) = x^\alpha \quad \text{and} \quad k(x) = 1 , \quad
\text{for } x>0 ,
\ee
where $\alpha \in \mbox{} ]0, 1[$ is a constant.
In this case, the function $\Theta$ defined by (\ref{Theta})
is given by
\be
\Theta (x) = x^\alpha - (r-b) x ,
\ee
and all of the conditions in Assumption~\ref{A5} hold true.
Furthermore,
\be
R_\Theta (x) = \frac{2}{\sigma^2 (\alpha-m) (n-\alpha)}
x^\alpha - x ,
\ee
which implies that
\be
\lim _{x \downarrow 0} \frac{R'_\Theta (x)}{\psi' (x)} =
\lim _{x \downarrow 0} \biggl( \frac{2\alpha}
{\sigma^2 n (\alpha-m) (n-\alpha)} x^{\alpha-n} -
\frac{1}{n} x^{1-n} \biggr)
= \infty
\ee
because $m < 0 < \alpha < 1 < n$.
In view of Lemmas~\ref{lem-z*} and~\ref{lem:FB}, we
can see that
\be
\xlb < \xub = c^\star = \infty .
\ee
Therefore, a $\beta$-$\gamma$ strategy is optimal
(Case~I of Theorem~\ref{thm:SolPb-lem}) for all
$c>0$.
\mbox{} \hfill $\Box$ }
\end{ex}

\begin{ex} \label{ex6} {\rm
Suppose that  $r+b-\sigma^2>0 \Leftrightarrow m<-1$
and $b>r$.
Also, consider the functions 
\begin{align*}
h(x) = \begin{cases} - \alpha x , & \text{if } x \in \mbox{} ]0,1[
, \\ - \alpha , \quad &\text{if } x \geq 1 , \end{cases}
\quad \text{and} \quad k(x) = \begin{cases} 3-2x , & \text{if }
x \in \mbox{} ]0,1[ , \\ x^{-2} , & \text{if } x \geq 1 , \end{cases}
\end{align*}
for some constant $\alpha \in \mbox{} \bigl] {-\infty} , 3(b-r)
\bigr[$.
In this case,
\be
\Theta (x) = \begin{cases} (r-2b-\sigma^2) x^2 + (3b-3r-\alpha)
x , & \text{if } x \in \mbox{} ]0,1[ , \\ (r+b-\sigma^2) x^{-1} -
3r - \alpha , & \text{if } x \geq 1, \end{cases}
\ee
and all of the conditions in Assumption~\ref{A5} hold true.
In view of the assumption that $m < -1$ and the identity
in (\ref{Q-Theta-x>>0}), we can see that $\lim
_{x \uparrow \infty} Q_\Theta (x) = \infty$, which implies
that $\xlb < \infty$ thanks to Lemma~\ref{lem-z*}.(ii).
Furthermore,
\be
R_\Theta (x) = R_h (x) - \int _0^x k(s) \, \di s =
\begin{cases} x^2 + \bigl( \frac{\alpha}{b-r} - 3\bigr) x -
\frac{2\alpha}{\sigma^2 (n-m) n (1-n)} x^n , &\text{if } x
\in \mbox{} ]0,1[ , \\ - \frac{\alpha}{r} - 3 + x^{-1} -
\frac{2\alpha}{\sigma^2 (n-m) m (1-m)} x^m , & \text{if }
x \geq 1 , \end{cases}
\ee
which implies that
\be
\lim _{x \downarrow 0} \frac{R'_\Theta (x)}{\psi' (x)} =
- \frac{2\alpha}{\sigma^2 (n-m) n (1-n)} \in \mbox{} \Bigl]
-\frac{6 (b-r)}{\sigma^2 (n-m) n (1-n)} , \infty \Bigr[ .
\ee
In view of Lemmas~\ref{lem-z*}, \ref{lem:FB}, \ref{lem:beta-0}
and~\ref{lem:oo-0}, we can see that, if $\alpha \leq 0$,
then
\be
\xub = \infty , \quad \text{and} \quad c^\star \in \mbox{}
]0, \infty[ ,
\ee
while, if $\alpha \in \mbox{} \bigl] 0 , 3(b-r) \bigr[$, then
\be
\xub < \infty \quad \text{and} \quad 0 < c^\star < c^\circ
= 3 + \frac{\alpha}{r} .
\ee
If $\alpha \leq 0$ and $c \in \mbox{} ]0,
c^\star[$, then a $\beta$-$\gamma$ strategy
is optimal (Case~I of Theorem~\ref{thm:SolPb-lem}), while,
if $\alpha \leq 0$ and $c \geq c^\star$, then no intervention
at all is optimal (Case~III of Theorem~\ref{thm:SolPb-lem}).
On the other hand, if $\alpha \in \mbox{} \bigl] 0 , 3(b-r)
\bigr[$, then any of the Cases I, II or III of
Theorem~\ref{thm:SolPb-lem} arises depending on whether
$0<c<c^\star$, $c^\star \leq c < 3 + \frac{\alpha}{r}$
or $c \geq 3 + \frac{\alpha}{r}$ is the case, respectively.
\mbox{} \hfill $\Box$ }
\end{ex}

\begin{ex} \label{ex7} {\rm
Suppose that $b = r > \frac{1}{2} \sigma^2$, which implies
that $m < -1$ and $n=1$.
Also, consider the functions 
\be
h(x) = \begin{cases} a x^\alpha , & \text{if } x \in \mbox{}
]0,1[ , \\ a , & \text{if } x \geq 1 , \end{cases} \quad
\text{and} \quad k(x) = \begin{cases} 4-2x , & \text{if } x
\in \mbox{} ]0,1[ , \\ x^{-2} + 1, & \text{if } x \geq 1 ,
\end{cases}
\ee
for some constants $\alpha \in \mbox ]1, 2]$ and
$a \in \mbox{} \bigl] 0, \frac{3}{2} (2b + \sigma^2) (\alpha-1)
\vee 3r \bigr[$. 
In this case, all of the conditions in Assumption~\ref{A5}
hold true,
\begin{gather}
\lim _{x \uparrow \infty} \psi ^{-1} (x) \int_0^\infty k(s) \di s
= 1, \quad \Theta (x) = \begin{cases} (r-2b-\sigma^2) x^2
+ a x^\alpha , & \text{if } x \in \mbox{} ]0,1[, \\ (r+b-\sigma^2)
x^{-1} - 3r + a , & \text{if } x \geq 1, \end{cases} 
\nonumber \\
\text{and} \quad
R_\Theta (x) = \begin{cases} x^2+\frac{2a}
{\sigma^2 (n-\alpha) (\alpha-m)} x^\alpha - \bigl( 3 +
\frac{2a \alpha}{\sigma^2 (n-m) n (n-\alpha)} \bigr) x ,
& \text{if } x \in \mbox{} ]0,1[ , \\ \frac{a}{r} - 3 + x^{-1} +
\frac{2a \alpha}{\sigma^2 (n-m) m (\alpha-m)} x^m ,
& \text{if } x \geq 1 . \end{cases} \nonumber
\end{gather}
Furthermore,
\be
\lim _{x \downarrow 0} \frac{R'_\Theta (x)}{\psi' (x)} =
-3 - \frac{2a}{(2b + \sigma^2) (\alpha-1)} < 0 .
\ee
In view of Lemmas~\ref{lem-z*}, \ref{lem:FB}, \ref{lem:beta-0}
and~\ref{lem:oo-0}, we can see that 
\be
\xlb < \xub < \infty \quad \text{and} \quad 0 < c^\star < c^\circ
= 3 - \frac{a}{r} .
\ee
Any of the Cases I, II or IV of Theorem~\ref{thm:SolPb-lem}
may arise, depending on whether $0<c<c^\star$, $c^\star \leq
c < 3 - \frac{a}{r}$ or $c \geq 3 - \frac{a}{r}$ is the case,
respectively.
\mbox{} \hfill $\Box$ }
\end{ex}

The next example shows that (\ref{R(0)-nat}) in
Example~\ref{lem:R_F} is not necessarily true if 0 is an
entrance boundary point.
Furthermore, it shows that $\beta$-0 strategies would
be an indispensable part of the optimal tactics if we
allowed for switching off the system and 0 were an
entrance boundary point (see Remark~\ref{rem:nat/entr}
at the end of the previous section).

\begin{ex} \label{ex8} {\rm
Suppose that $X$ is the mean-reverting square-root process
that is given by
\be
\di X_t = \alpha (2 - X_t) \, \di t +  \sqrt{2 \alpha X_t}
\, \di W_t , \quad X_0 = x>0 ,
\ee
for some constant $\alpha > 0$.
Also, suppose that
\be
r(x) = \alpha , \quad h(x) = \begin{cases} e^x -1 , &
\text{if } x \in \mbox{} ]0, 1[ , \\ e - e^{\gamma+3}
- 1 + e^{\gamma x +3} , & \text{if } x \geq 1 , \end{cases}
\quad \text{and} \quad k(x) = \kappa , \quad
\text{for } x>0 ,
\ee
for some constants $\gamma < 0$ and
$\kappa \in \mbox{} \bigl] 0, \frac{1}{2 \alpha} \bigr[$.
In this case,
\be
\varphi (x) = \frac{1}{x} , \quad \psi (x) = \frac{e^x - 1}{x}
\quad \text{and} \quad p'(x) = \frac{1}{x^2} e^{x-1} .
\ee
In particular, 0 is an entrance boundary point.
The function $\Theta$ defined by (\ref{Theta})
is given by
\be
\Theta (x) =  \begin{cases} 2 \alpha \kappa - 1 -
2 \alpha \kappa x + e^x , & \text{if } x \in \mbox{}
]0, 1[ , \\ 2 \alpha \kappa + e - e^{\gamma+3}
- 1 - 2 \alpha \kappa x + e^{\gamma x +3} , & \text{if } x \geq 1
, \end{cases}
\ee
all of the conditions in Assumption~\ref{A5} hold true,
\be
\lim _{x \downarrow 0} R_h (x) = \frac{1}{\alpha} \biggl(
1 + \frac{\gamma}{1-\gamma} e^{\gamma+2} \biggr)
=: \frac{1}{\alpha} f(\gamma) \quad \text{and} \quad
\lim _{x \downarrow 0} \frac{R'_\Theta (x)}{\psi' (x)} =
\frac{1}{\alpha} f(\gamma) - 2 \kappa .
\ee
The function $f$ is strictly decreasing in
the interval $\bigl] {-\infty} , (1-\sqrt{5}) / 2 \bigr[$,
strictly increasing in the interval $\bigl] (1-\sqrt{5})
/ 2 , 0 \bigr[$,
\be
\lim _{\gamma \downarrow -\infty} f(\gamma)
= 1, \quad f \biggl( \frac{1-\sqrt{5}}{2} \biggr)
= 1 - \frac{1}{2} \bigl( 3 - \sqrt{5} \bigr)
e^{(5 - \sqrt{5}) / 2} < 0 \quad \text{and} \quad
f(0) = 1 .
\ee
Therefore, there exist constants $(1-\sqrt{5}) / 2 <
\gamma _1 < \gamma_2 < 0$ such that $f(\gamma)
< 0$ for all $\gamma \in \bigl[ (1-\sqrt{5})
/ 2 , \gamma_1 \bigr[$ and $f(\gamma) \in \mbox{}
]0, 2 \alpha \kappa[$ for all $\gamma \in \mbox{}
]\gamma _1, \gamma_2[$.
In view of these observations, we can see that
\ben
\lim _{x \downarrow 0} R_h (x) \neq 0 =
\lim _{x \downarrow 0} \frac{h(x)}{r(x)} \quad \text{ for all }
\gamma \in \bigl[ (1-\sqrt{5}) / 2 , \gamma_1 \bigr[
\mbox{} \setminus \{ \gamma_2 \} , \label{Rh(0)}
\een
which shows that (\ref{R(0)-nat}) in Lemma~\ref{lem:R_F}
is not in general true if 0 is an entrance boundary point.
On the other hand, Lemma~\ref{lem-z*} implies that,
if $\gamma \in \mbox{} \bigl[ (1-\sqrt{5}) / 2 , \gamma_1
\bigr[$, then $0 < \xlb < \xub < \infty$ and we are in the
context of Lemma~\ref{lem:beta-0} with $c^\circ = \infty$.
In this context, (\ref{w-sol-0}) yields the expression
\be
w(0) = \lim _{x \downarrow 0} w(x) = \frac{1}{\alpha}
f(\gamma) - \frac{R_\Theta' (\beta^\circ)}{\psi' (\beta^\circ)} .
\ee
In view of (\ref{z_*}) in Lemma~\ref{lem-z*}, (\ref{beta^o})
in Lemma~\ref{lem:beta-0}, Remark~\ref{rem:nat/entr} and
the analysis thus far, we can see the following:

\noindent
(a) If $\gamma \in \mbox{} ]\gamma_1, 0[$, then
$w(0) > 0 = h(0) / r(0)$ and a $\beta$-0 strategy would be
strictly sub-optimal.

\noindent
(b) If $\gamma \in \mbox{} \bigl] (1-\sqrt{5}) / 2 ,
\gamma_1 \bigr[$, then $w(0) < 0 = h(0) / r(0)$ for
all $c$ sufficiently large, in which case, a
$\beta$-0 strategy would be optimal.
\mbox{} \hfill $\Box$ }
\end{ex}

Our final example shows that the conditions in
(\ref{xlb<oo-SC}) are only sufficient for the point $\xlb$
introduced in part~(i) of Lemma~\ref{lem-z*} to
be finite.

\begin{ex} \label{ex9} {\rm
Suppose that $X$ is the geometric Brownian motion
that is given by (\ref{gBm}) with $b = \frac{1}{4}$ and
$\sigma = \frac{1}{\sqrt{2}}$.
Also, suppose that $r=1$, so that $m = -2$, $n=2$
and $C \equiv n-m = 4$.
The functions defined by
\be
k(x) = \begin{cases} 6-5x , & \text{if } x \leq 1 , \\
x^{-5} , & \text{if } x > 1 , \end{cases}
\quad \text{and} \quad
h(x) = \begin{cases} 7x , & \text{if } x \leq 1 , \\
6 + x^{-4} , & \text{if } x > 1 , \end{cases}
\ee
are such that
\be
\Theta (x) = \begin{cases} \frac{5}{2} x , & \text{if }
x \leq 1 , \\ \frac{9}{4} + \frac{1}{4} x^{-4} , & \text{if }
x > 1 , \end{cases}
\ee
and the function $Q_\Theta$ defined by
(\ref{Q-Theta}) satisfies $\lim _{x \uparrow \infty}
Q_\Theta (x) = - \frac{1}{6}$.
In this case, the necessary and sufficient condition
of Lemma~\ref{lem-z*}.(ii) implies that $\xlb = \infty$.
On the other hand, the functions defined by
\be
k(x) = \begin{cases} 6-5x , & \text{if } x \leq 1 , \\
x^{-5} , & \text{if } x > 1 , \end{cases}
\quad \text{and} \quad
h(x) = \begin{cases} 5x , & \text{if } x \leq 1 , \\
4 + x^{-4} , & \text{if } x > 1 , \end{cases}
\ee
are such that
\be
\Theta (x) = \begin{cases} \frac{1}{2} x , & \text{if }
x \leq 1 , \\ \frac{1}{4} + \frac{1}{4} x^{-4} , & \text{if }
x > 1 , \end{cases}
\ee
and the function $Q_\Theta$ defined by
(\ref{Q-Theta}) satisfies $\lim _{x \uparrow \infty}
Q_\Theta (x) = \frac{1}{6}$.
In this case, the necessary and sufficient condition
of Lemma~\ref{lem-z*}.(ii) implies that $\xlb < \infty$.
\mbox{} \hfill $\Box$ }
\end{ex}

%%%%%%%%%%%%%%%%%%%%%%%%%%%%%%%%%%


\begin{thebibliography}{20}
%
\bibitem{A041}
{\sc  L.\,H.\,R.\,Alvarez} (2004),
A class of solvable impulse control problems,
{\em Applied Mathematics and Optimization\/}, {\bf 49}, 
pp.\,265--295.
%
\bibitem{A042}
{\sc  L.\,H.\,R.\,Alvarez} (2004),
Stochastic forest stand value and optimal timber harvesting,
{\em SIAM Journal of Control and Optimization\/}, {\bf 42}, 
pp.\,1972--1993.
%
\bibitem{AH22}
{\sc  L.\,H.\,R.\,Alvarez and E.\,A.\,Hening} (2022),
Optimal sustainable harvesting of populations
in random environments,
{\em Stochastic Processes and their Applications\/},
vol.\,{\bf 150}, pp.\,678--698.
%
\bibitem{AK07}
{\sc  L.\,H.\,R.\,Alvarez and E.\,Koskela} (2007),
The forest rotation problem with stochastic harvest and
amenity value,
{\em Natural Resource Modeling\/}, {\bf 20},  pp.\, 477--509.
%
\bibitem{AL08}
{\sc  L.\,H.\,R.\,Alvarez and J.\,Lempa} (2008),
On the optimal stochastic impulse control of linear diffusions,
{\em SIAM Journal of Control and Optimization\/}, {\bf 47},
pp.\,703--732. 
%
\bibitem{BPS04}
{\sc A.\,Bar-Ilan, D.\,Perry and W.\,Stadje} (2004),
A generalized impulse control model of cash management,
{\em Journal of Economic Dynamics \& Control\/}, vol.\,{\bf 28},
pp.\,1013--1033.
%
\bibitem{BSZ02}
{\sc A.\,Bar-Ilan, A.\,Sulem and A.\,Zanello} (2002),
Time-to-build and capacity choice,
{\em Journal of Economic Dynamics \& Control\/}, vol.\,{\bf 26},
pp.\,69--98.
%
\bibitem{BL84}
{\sc A.\,Bensoussan and J.\,L.\,Lions} (1984),
{\em Impulse Control and Quasivariational Inequalities\/}
Gauthier-Villars, Montrouge.
%
\bibitem{BP00}
{\sc T.\,Bielecki and S.\,Pliska} (2000),
Risk sensitive asset management with transaction costs,
{\em Finance \& Stochastics\/}, vol.\,{\bf 4}, pp.\,1--33.
%
\bibitem{BS}
{\sc  A.\,N.\,Borodin and P.\,Salminen} (2015),
{\em Handbook of Brownian Motion - Facts and Formulae\/},
Birkh{\"a}user, pp.\,37--38.
%
\bibitem{C00}
{\sc A.\,Cadenillas} (2000),
Consumption-investment problems with transaction costs:
survey and open problems,
{\em Mathematical Methods of Operations Research\/}, vol.\,{\bf 51},
pp.\,43--68.
%
\bibitem{CSZ07}
{\sc A.\,Cadenillas, S.\,Sarkar and F.\,Zapatero} (2007),
Optimal dividend policy with mean-reverting cash reservoir,
{\em Mathematical Finance\/}, vol.\,{\bf 17}, pp.\,81--109.
%
\bibitem{C14}
{\sc S.\,Christensen} (2014),
On the solution of general impulse control problems
using superharmonic functions,
{\em Stochastic Processes and their Applications\/},
vol.\,{\bf 124}, pp.\,709--729.
%
\bibitem{CS23}
{\sc S.\,Christensen and C.\,Strauch} (2023),
Nonparametric learning for impulse control
problems--exploration vs.\ exploitation,
{\em The Annals of Applied Probability\/}, vol.\,{\bf 33},
pp.\,1569--1587.
%
\bibitem{D93}
{\sc M.\,H.\,A.\,Davis} (1993),
{\em Markov Models and Optimization\/},
Chapman \& Hall.
%
\bibitem{DGW10}
{\sc M.\,H.\,A.\,Davis, X.\,Guo and G.\,Wu} (2009),
Impulse control of multidimensional jump
diffusions,
{\em SIAM Journal on Control and Optimization\/},
vol.\,{\bf 48}, pp.\,5276--5293.
%
\bibitem{DHH10}
{\sc B.\,Djehiche, S.\,Hamad\`{e}ne and I.\,Hdhiri} (2010),
Stochastic impulse control of non-Markovian processes,
{\em Applied Mathematics and Optimization\/}, vol.\,{\bf 61},
1--26.
%
\bibitem{E08}
{\sc M.\,Egami} (2008),
A direct solution method for stochastic impulse control
problems of one-dimensional diffusions,
{\em SIAM Journal on Control and Optimization\/},
vol.\,{\bf 47}, pp.\,1191--1218.
%
\bibitem{F03}
{\sc N.\,C.\,Framstad} (2003),
Optimal harvesting of a jump diffusion population and the
effect of jump uncertainty,
{\em SIAM Journal on Control and Optimization\/}, {\bf 42},
pp.\,1451--1465.
%
\bibitem{HST83}
{\sc J.\,M.\,Harrison, T.\,M.\,Sellke and A.\,Taylor} (1983),
Impulse control of Brownian motion,
{\em Mathematics of Operations Research\/}, vol.\,{\bf  8},
pp.\,454--466. 
%
\bibitem{HT83}
{\sc J.\,M.\,Harrison and M.\,I.\,Taksar} (1983),
Instantaneous control of Brownian motion,
{\em Mathematics of Operations Research\/}, vol.\,{\bf 8},
pp.\,439--453.
%
\bibitem{HSZ15}
{\sc K.\,L.\,Helmes, R.\,H.\,Stockbridge and C.\,Zhu} (2015),
A measure approach for continuous inventory models:
discounted cost criterion,
{\em SIAM Journal on Control and Optimization\/},
vol.\,{\bf 53}, pp.\,2100--2140.
%
\bibitem{HSZ24}
{\sc K.\,L.\,Helmes, R.\,H.\,Stockbridge and C.\,Zhu} (2024),
On the modelling of impulse control with random
effects for continuous Markov processes,
{\em SIAM Journal on Control and Optimization\/},
vol.\,{\bf 62}, pp.\,699--723.
%
\bibitem{K98}
{\sc R.\,Korn} (1998),
Portfolio optimization with strictly positive transaction
costs and impulse control,
{\em Finance \& Stochastics\/}, vol.\,{\bf 2}, pp.\,85--114.
%
\bibitem{K99}
{\sc R.\,Korn} (1999),
Some applications of impulse control in mathematical finance,
{\em Mathematical Methods of Operations Research\/}, vol.\,{\bf 50},
pp.\,493--518.
%
\bibitem{K80}
{\sc N.\,V.\,Krylov} (1980),
{\em Controlled Diffusion Processes}, Springer.
%
\bibitem{LZ}
{\sc D.\,Lamberton and M.\,Zervos} (2013),
On the optimal stopping of a one-dimensional diffusion,
{\em Electronic Journal of Probability\/}, {\bf 18} pp.\,1--49.  
%
\bibitem{LM84}
{\sc J.\,P.\,Lepeltier and B.\,Marchal} (1984),
General theory of Markov impulse control
{\em SIAM Journal on Control and Optimization\/}, vol.\,{\bf 22},
pp.\,645--665.
%
\bibitem{LO01}
{\sc E.\,Lungu and B.\,{\O}ksendal} (2001),
Optimal harvesting from interacting populations in a stochastic
environment, {\em Bernoulli\/}, {\bf 7}, pp.\,527--539.
%
\bibitem{LVMP07}
{\sc V.\,LyVath, M.\,Mnif and H.\,Pham} (2007),
A model of optimal portfolio selection under liquidity risk
and price impact, {\em Finance \& Stochastics\/},
vol.\,{\bf 11}, pp.\,51--90.
%
\bibitem{MR17}
{\sc J.\,L.\,Menaldi and M.\,Robin} (2017),
On some impulse control problems with constraint,
{\em SIAM Journal on Control and Optimization\/},
vol.\,{\bf 55}, pp.\,3204--3225.
%
\bibitem{MO98}
{\sc G.\,Mundaca and B.\,{\O}ksendal} (1998), 
Optimal stochastic intervention control with application to the
exchange rate, {\em Journal of Mathematical Economics\/},
{\bf 29}, pp.\,225--243.
%
\bibitem{OT06}
{\sc M.\,Ohnishi and M.\,Tsujimura} (2006),
An impulse control of a geometric Brownian motion with quadratic
costs, {\em European Journal of Operational Research\/},
vol.\,{\bf 168}, pp.\,311--321.
%
\bibitem{OS07}
{\sc B.\,{\O}ksendal and A.\,Sulem} (2007),
{\em Applied Stochastic Control of Jump Diffusions\/},
2nd edition, Springer.
%
\bibitem{PS17}
{\sc J.\,Palczewski and \L.\,Stettner} (2017),
Impulse control maximizing average cost per unit
time: a nonuniformly ergodic case,
 {\em SIAM Journal on Control and Optimization\/},
vol.\,{\bf 55}, pp.\,936--960.
%
\bibitem{P84}
{\sc B.\,Perthame} (1984), Continuous and impulsive control of
diffusion processes in $\bbr ^N$, {\em Nonlinear Analysis\/},
{\bf 8}, pp.\,1227--1239. 
%
\bibitem{P09}
{\sc H.\,Pham} (2009),
{\em Continuous-time Stochastic Control and Optimization with
Financial Applications\/}, Springer.
%
\bibitem{RY}
{\sc D.\,Revuz and M.\,Yor} (1999),
{\em Continuous Martingales and Brownian Motion\/},
3rd edition, Springer.
%
\bibitem{R77}
{\sc S.\,Richard} (1977),
Optimal impulse control of a diffusion process with
both fixed and proportional costs,
{\em SIAM Journal on Control and Optimization\/},
vol.\,{\bf 15}, pp.\,79--91.
%
\bibitem{SSZ11}
{\sc Q.\,Song, R.\,H.\,Stockbridge and C.\,Zhu} (2011),
On optimal harvesting problems in random environments,
{\em SIAM Journal on Control and Optimization\/}, {\bf 49},
pp.\,859--889.
%
\bibitem{S83}
{\sc \L.\,Stettner} (1983),
On impulsive control with long run average cost criterion,
{\em Studia Mathematica\/}, vol.\,{\bf 76}, pp.\,279--298.
%
\end{thebibliography}
\end{document}